\documentclass[a4paper,reqno]{amsart}
\usepackage[british]{babel}
\usepackage[T1]{fontenc}
\usepackage{textcomp}
\usepackage{mathtools}  
\usepackage{amsthm}
\usepackage{amssymb}
\usepackage{mathptmx} 
\usepackage{microtype}
\usepackage{paralist}
\usepackage{nicefrac}
\usepackage{xcolor}
\usepackage{tikz}
  \usetikzlibrary{arrows,positioning}
  \definecolor{UGentblauw}{cmyk}{1,.8,.3,.05}
  \definecolor{UGentgeel}{cmyk}{0,.3,1,0}
  \tikzstyle{root}=[rectangle,draw=blue!90]
  \tikzstyle{nonterminal}=[rectangle,rounded corners,fill=blue!15,draw=blue!15]
  \tikzstyle{terminal}=[rectangle]
  \tikzstyle{cut}=[thick,dotted,draw=green!50!black]
  \tikzstyle{local}=[color=green!50!black,text=green!25!black]
  \tikzstyle{simplexbackground}=[UGentblauw,fill=UGentblauw!5,text=black]
  \tikzstyle{simplexborder}=[UGentblauw!60]
  \tikzstyle{bluebackground}=[very thin,draw=blue,fill=blue,fill opacity=0.6]
  \tikzstyle{yellowbackground}=[very thin,draw=darkgray,fill=yellow,fill opacity=1]
  \tikzstyle{redbackground}=[very thin,draw=red,fill=red,fill opacity=0.6]
  \tikzstyle{statesbackground}=[draw=black!40!UGentblauw!50,fill=UGentblauw!10]
  \tikzstyle{closedbackground}=[draw=black!40!UGentblauw!50,fill=white]
  \tikzstyle{communbackground}=[draw=black!40!UGentgeel!75,fill=UGentgeel!30]
  \tikzstyle{cyclic}=[style=help lines]
  \newcommand*{\simplexpath}{(0,1,0) -- (1,0,0) -- (0,0,1) -- cycle}
\usepackage[longnamesfirst,numbers,sort&compress]{natbib}
  \newcommand{\citetopt}[2][]{\citeauthor{#2} \citep[#1]{#2}} 
\usepackage[bookmarks,bookmarksnumbered]{hyperref}  
\theoremstyle{plain}
  \newtheorem{definition}{Definition}[section]
  \newtheorem{theorem}{Theorem}[section]
  \newtheorem{proposition}[theorem]{Proposition}
  \newtheorem{corollary}[theorem]{Corollary}
  \newtheorem{lemma}[theorem]{Lemma}
\theoremstyle{remark}
  \newtheorem{example}{Example}[section]
\DeclareRobustCommand{\newconcept}[1]{\emph{#1}}
\DeclareRobustCommand{\pties}{\mathcal{Y}} 
\DeclareRobustCommand{\pty}{y} 
\DeclareRobustCommand{\ex}{E} 
\DeclareRobustCommand{\lex}{\underline{\ex}\vphantom{\ex}} 
\DeclareRobustCommand{\uex}{\overline{\ex}\vphantom{\ex}} 
\DeclareRobustCommand{\trans}{\mathrm{T}} 
\DeclareRobustCommand{\utrans}{\overline{\trans}\vphantom{\trans}}  
\DeclareRobustCommand{\ltrans}{\underline{\trans}\vphantom{\trans}}  
\DeclareRobustCommand{\ttrans}{\mathbb{T}}  
\DeclareRobustCommand{\uttrans}{\overline{\ttrans}\vphantom{\ttrans}}  
\DeclareRobustCommand{\transmats}{\mathcal{T}}  
\DeclareRobustCommand{\relty}{r}  
\DeclareRobustCommand{\lrelty}{\underline{\relty}}  
\DeclareRobustCommand{\urelty}{\overline{\relty}}  

\DeclareRobustCommand{\utp}[3][]{\overline{P}\vphantom{P}^{#1}_{\mspace{-8mu}#2#3}} 
\DeclareRobustCommand{\transmat}{T} 
\DeclareRobustCommand{\utransmat}{\overline{\transmat}} 
\DeclareRobustCommand{\ltransmat}{\underline{\transmat}} 
\DeclareRobustCommand{\dismat}{m} 
\DeclareRobustCommand{\udismat}{\overline{\dismat}\vphantom{\dismat}} 
\DeclareRobustCommand{\ldismat}{\underline{\dismat}\vphantom{\dismat}} 

\DeclareRobustCommand{\states}{{\mathcal{X}}} 
\DeclareRobustCommand{\maxregstates}{{\mathcal{R}}} 
\DeclareRobustCommand{\allgambles}{\mathcal{L}} 
\DeclareRobustCommand{\margmass}{\mathcal{M}} 
\DeclareRobustCommand{\mass}{\mathcal{P}} 
\DeclareRobustCommand{\condmass}{\mathcal{Q}} 
\DeclareRobustCommand{\expects}{\mathcal{E}} 

\DeclareRobustCommand{\simplex}{\Sigma}
\DeclareRobustCommand{\domain}{\mathcal{K}}
\DeclareMathOperator{\ext}{ext}

\DeclareRobustCommand{\marg}{m} 
\DeclareRobustCommand{\umarg}{\overline{\marg}\vphantom{\marg}} 
\DeclareRobustCommand{\lmarg}{\underline{\marg}\vphantom{\marg}} 
\DeclareRobustCommand{\cond}{q} 
\DeclareRobustCommand{\ucond}{\overline{\cond}\vphantom{\cond}} 
\DeclareRobustCommand{\lcond}{\underline{\cond}\vphantom{\cond}} 

\DeclareRobustCommand{\commun}[2]{#1\leftrightsquigarrow#2}
\DeclareRobustCommand{\access}[3][]{\smash[t]{#2\stackrel{#1}\rightsquigarrow#3}}

\DeclareRobustCommand{\uaccess}[3][]{\smash[t]{#2\stackrel{#1}\rightarrow#3}}

\DeclareRobustCommand{\nsteps}[2]{N_{#1#2}}
\DeclareRobustCommand{\period}[1]{d_{#1}}
\DeclareRobustCommand{\steps}[2]{t_{#1#2}}
\DeclareRobustCommand{\dif}{\,\mathrm{d}} 

\DeclareRobustCommand{\then}{\Rightarrow}

\DeclareRobustCommand{\ftuple}[3][1]{({#2}_{#1:#3})} 
\DeclareRobustCommand{\ntuple}[3][1]{{#2}_{#1:#3}} 
\DeclareRobustCommand{\vtuple}[3][1]{{#2}_{#1:#3}} 

\DeclareRobustCommand{\set}[2]{\left\{#1\colon#2\right\}} 
\DeclareRobustCommand{\inlineset}[2]{\{#1\colon#2\}} 
\DeclareRobustCommand{\asa}{\Leftrightarrow} 
\DeclareRobustCommand{\reals}{\mathbb{R}} 

\DeclareRobustCommand{\naturals}{\mathbb{N}} 
\DeclareRobustCommand{\ind}[1]{I_{#1}} 
\DeclareRobustCommand{\cg}[1][\states]{\ind{#1}} 
\DeclareRobustCommand{\init}{\square} 
\DeclareRobustCommand{\eqdef}{\coloneqq} 
\DeclareRobustCommand{\defeq}{\eqqcolon} 

\DeclarePairedDelimiter{\supnorm}{\lVert}{\rVert_\infty}
\DeclarePairedDelimiter{\abs}{\lvert}{\rvert}
\DeclarePairedDelimiter{\card}{\lvert}{\rvert}
\DeclareMathOperator{\id}{id} 
\begin{document}
\title{Imprecise Markov chains and their limit behaviour}

\author{Gert de Cooman}
\author{Filip Hermans}
\author{Erik Quaeghebeur}
\address{SYSTeMS Research Group, Ghent University, Technologiepark--Zwijnaarde 914, 9052 Zwijnaarde, Belgium}
\email{\{gert.decooman,filip.hermans,erik.quaeghebeur\}@UGent.be}




\begin{abstract}
  When the initial and transition probabilities of a finite Markov chain in discrete time are not well known, we should perform a sensitivity analysis.
  This can be done by considering as basic uncertainty models the so-called \newconcept{credal sets} that these probabilities are known or believed to belong to, and by allowing the probabilities to vary over such sets.
  This leads to the definition of an \newconcept{imprecise Markov chain}.
  We show that the time evolution of such a system can be studied very efficiently using so-called \newconcept{lower} and  \newconcept{upper expectations}, which are equivalent mathematical representations of credal sets.
  We also study how the inferred credal set about the state at time $n$ evolves as $n\to\infty$: under quite unrestrictive conditions, it converges to a uniquely invariant credal set, regardless of the credal set given for the initial state.
  This leads to a non-trivial generalisation of the classical Perron--Frobenius Theorem to imprecise Markov chains.
\end{abstract}

\maketitle

\section{Introduction}
One convenient way to model uncertain dynamical systems is to describe them as Markov chains. These have been studied in great detail, and their properties are well known. However, in many practical situations, it remains a challenge to accurately identify the transition probabilities in the Markov chain: the available information about physical systems is often imprecise and uncertain. Describing a real-life dynamical system as a Markov chain will therefore often involve unwarranted precision, and may lead to conclusions not supported by the available information.
\par
For this reason, it seems quite useful to perform probabilistic robustness studies, or sensitivity analyses, for Markov chains. This is especially relevant in decision-making applications. Many researchers in Markov Chain Decision Making \citep{white1994,harmanec2002,nilim2005,itoh2007}---inspired by \citeauthor{satia1973}'s~\citeyearpar{satia1973} original work---have paid attention to this issue of `imprecision' in Markov chains.
\par
Work on the more mathematical aspects of modelling such imprecision in Markov chains was initiated in the early 1980s by \citeauthor{hartfiel1994} (see \cite{hartfiel1991,hartfiel1994,hartfiel1998}), under the name `Markov set-chains'.  \Citeauthor{hartfiel1991}'s work seems to have been unknown to \citet{kozine2002}, who approached the subject from a different angle. 
Armed with linear programming techniques, these authors performed an experimental study of the limit behaviour of Markov chains with uncertain transition probabilities. 
More recently, \citet{skulj2006,skulj2007} has also contributed to a formal study of the time evolution and limit behaviour of such systems. Markov set-chains can also be seen as special cases of so-called \emph{credal networks} under strong independence \cite{cozman2000,cozman2005}.  
\par
All these approaches use \newconcept{sets of probabilities} to deal with the imprecision in the transition probabilities. When these probabilities are not well known, they are assumed to belong to certain sets, and robustness analyses are performed by allowing the transition probabilities to vary over such sets. This should be contrasted with more common ways of performing a sensitivity analysis: looking at small deviations from a reference model and evaluating derivatives of important variables in this reference point. 
\par
As we shall see, the sets of probabilities approach leads to a number of computational difficulties. But we will show that they can be overcome by tackling the problem from another angle, using lower and upper expectations, rather than sets of probabilities. Our new method also makes it fairly easy to formulate and prove convergence (or Perron--Frobenius-like) results for  Markov chains with uncertain transition probabilities that hold under weaker conditions than the ones found by \citet{hartfiel1991,hartfiel1998} and \citet{skulj2007}. We shall see that our condition for this convergence, which requires that the imprecise Markov chain should be \emph{regularly absorbing}, is implied by, and even strictly weaker than, both \citeauthor{hartfiel1998}'s \emph{product scrambling} and \citeauthor{skulj2007}'s \emph{regularity} conditions.
\par
In the rest of this Introduction, we give an overview of the theory of classical Markov chains and formulate the classical Perron--Frobenius theorem. Then, in Sections~\ref{sec:towards} and~\ref{sec:sensitivity-analysis}, we introduce imprecise Markov chains and generalise many aspects of the classical theory. In Section~\ref{sec:accessibility}, we briefly discuss accessibility relations, which allows us to give a nice interpretation to a number of conditions that will turn out to be sufficient for a Perron--Frobenius-like convergence result. In Section~\ref{sec:convergence}, we generalise the classical Perron--Frobenius theorem, and explore the relation of our generalisation with previous work in the literature. We discuss a number of theoretical and numerical examples in Section~\ref{sec:examples}, and we give perspectives for further research in the Conclusions. Proofs of theorems and propositions have been relegated to an appendix.

\subsection{A short analysis of classical Markov chains}
\begin{figure}[ht]
  \centering\footnotesize
  \begin{tikzpicture}
    \tikzstyle{level 1}=[sibling distance=20em]
    \tikzstyle{level 2}=[sibling distance=10em]
    \tikzstyle{level 3}=[sibling distance=5em]
    \node[root] (root) {} [grow=down,level distance=8ex]
    child {node[nonterminal] (a) {$a\vphantom{)}$}
      child {node[nonterminal] (aa) {$(a,a)$}
        child {node[nonterminal] (aaa) {$(a,a,a)$}}
        child {node[nonterminal] (aab) {$(a,a,b)$}}
      }
      child {node[nonterminal] (ab) {$(a,b)$}
        child {node[nonterminal] (aba) {$(a,b,a)$}}
        child {node[nonterminal] (abb) {$(a,b,b)$}}
      }
    }
    child {node[nonterminal] (b) {$b\vphantom{)}$}
      child {node[nonterminal] (ba) {$(b,a)$}
        child {node[nonterminal] (baa) {$(b,a,a)$}}
        child {node[nonterminal] (bab) {$(b,a,b)$}}
      }
      child {node[nonterminal] (bb) {$(b,b)$}
        child {node[nonterminal] (bba) {$(b,b,a)$}}
        child {node[nonterminal] (bbb) {$(b,b,b)$}}
      }
    };
    \draw[cut] (b) -- +(1,0);
    \draw[cut] (b) -- (a) -- +(-2,0) node[left,local] {$\states^1$};
    \draw[cut] (bb) -- +(1,0);
    \draw[cut] (bb) -- (ba) -- (ab) -- (aa) -- +(-1,0) node[left,local] {$\states^2$};
  \end{tikzpicture}
  \caption{
    The event tree for the time evolution of system that can be in two states, $a$~and~$b$, and can change state at time instants $n=1,2$.
    Also depicted are the respective cuts~$\states^1$ and~$\states^2$ of\/~$\init$ where the states at times~$1$ and~$2$ are revealed.}
  \label{fig:markov-event}
\end{figure}
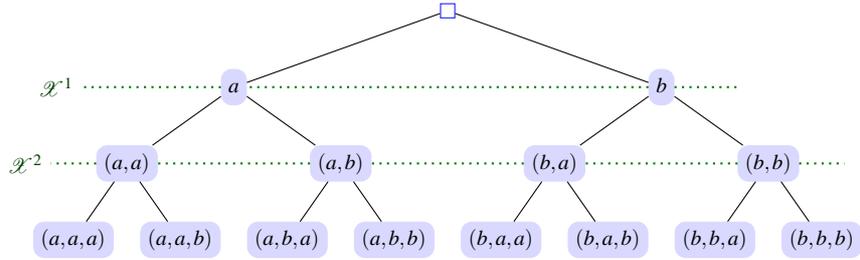

Consider a finite Markov chain in discrete time, where at consecutive times $n=1,2,3,\dots,N$, $N\in\naturals$ the \newconcept{state}~$X(n)$ of a system can assume any value in a finite set~$\states$. Here $\naturals$~denotes the set of non-zero natural numbers, and~$N$ is the time horizon. The time evolution of such a system can be modelled as if it traversed a path in a so-called \newconcept{event tree}; see \citet{shafer1996a}. An example of such a tree for $\states=\{a,b\}$ and $N=3$ is given in Figure~\ref{fig:markov-event}.
\par
The \newconcept{situations}, or nodes, of the tree have the form $\vtuple{x}{k}\eqdef(x_1,\ldots,x_k)\in\states^k$, $k=0,1,\dots,N$. For $k=0$ there is some abuse of notation as we let $\states^0\eqdef\{\init\}$, where~$\init$ is the so-called \newconcept{initial situation}, or root of the tree. In the cuts\footnote{A \newconcept{cut} $V$ of a situation $s$ is a collection of descendants $v$ of $s$ such that every path (from root to leaves) through $s$ goes through exactly one $v$ in $V$.} $\states^n$ of~$\init$, the value of the state $X(n)$ at time $n$ is revealed.
\par
\par
In a classical analysis, it is generally assumed that we have: (i) a probability distribution over the initial state $X(1)$, in the form of a probability mass function $m_1$ on $\states$; and (ii) for each situation $\vtuple{x}{n}$ that the system can be in at time $n$, a~probability distribution over the next state $X(n+1)$, in the form of a probability mass function $q(\cdot\vert\vtuple{x}{n})$ on $\states$. This means that in each non-terminal situation\footnote{A \newconcept{non-terminal} situation is a node of the tree that is not a   leaf.} $\vtuple{x}{n}$ of the event tree, we have a \emph{local} probability model telling us about the probabilities of each of its child nodes. This turns the event tree into a so-called \newconcept{probability tree}; see \citetopt[Chapter~3]{shafer1996a} and \citetopt[Section~1.9]{kemeny1976}.
\par
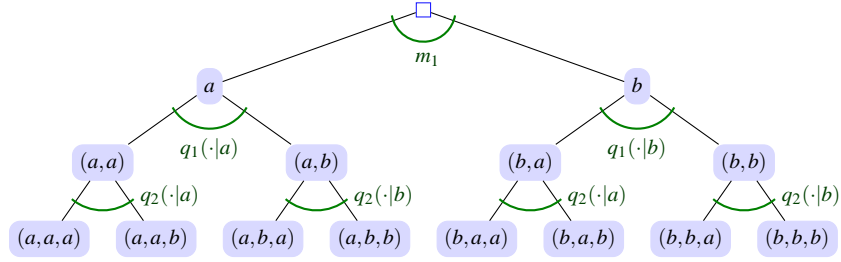
\begin{figure}[ht]
  \centering\footnotesize
  \begin{tikzpicture}
    \tikzstyle{level 1}=[sibling distance=20em]
    \tikzstyle{level 2}=[sibling distance=10em]
    \tikzstyle{level 3}=[sibling distance=5em]
    \node[root] (root) {} [grow=down,level distance=8ex]
    child {node[nonterminal] (a) {$a\vphantom{)}$}
      child {node[nonterminal] (aa) {$(a,a)$}
        child {node[nonterminal] (aaa) {$(a,a,a)$}}
        child {node[nonterminal] (aab) {$(a,a,b)$}}
      }
      child {node[nonterminal] (ab) {$(a,b)$}
        child {node[nonterminal] (aba) {$(a,b,a)$}}
        child {node[nonterminal] (abb) {$(a,b,b)$}}
      }
    }
    child {node[nonterminal] (b) {$b\vphantom{)}$}
      child {node[nonterminal] (ba) {$(b,a)$}
        child {node[nonterminal] (baa) {$(b,a,a)$}}
        child {node[nonterminal] (bab) {$(b,a,b)$}}
      }
      child {node[nonterminal] (bb) {$(b,b)$}
        child {node[nonterminal] (bba) {$(b,b,a)$}}
        child {node[nonterminal] (bbb) {$(b,b,b)$}}
      }
    };
    \draw[local,thick] (root) +(190:1.5em) arc (190:350:1.5em);
    \draw[local,thick] (b) +(210:2em) arc (210:330:2em);
    \draw[local,thick] (a) +(210:2em) arc (210:330:2em);
    \draw[local,thick] (bb) +(230:2.25em) arc (230:310:2.25em);
    \draw[local,thick] (ba) +(230:2.25em) arc (230:310:2.25em);
    \draw[local,thick] (ab) +(230:2.25em) arc (230:310:2.25em);
    \draw[local,thick] (aa) +(230:2.25em) arc (230:310:2.25em);
    \path (root) +(275:2.35em) node[local] {$m_1$};
    \path (a) +(270:2.95em) node[local] {$q_1(\cdot\vert a)$};
    \path (b) +(270:2.95em) node[local] {$q_1(\cdot\vert b)$};
    \path (aa) +(300:2.8em) node[local,above right] {$q_2(\cdot\vert a)$};
    \path (bb) +(300:2.8em) node[local,above right] {$q_2(\cdot\vert b)$};
    \path (ba) +(300:2.8em) node[local,above right] {$q_2(\cdot\vert a)$};
    \path (ab) +(300:2.8em) node[local,above right] {$q_2(\cdot\vert b)$};
  \end{tikzpicture}
  \caption{The probability tree for the time evolution of a~Markov chain that can be in two states, $a$ and $b$, and can change state at each time instant $n=1,2$.}
  \label{fig:markov-probability}
\end{figure}

The probability tree for a Markov chain is special, because the \newconcept{Markov Condition} states that when the system jumps from state $X(n)=x_n$ to a new state $X(n+1)$, where the system goes to will only depend on the state $X(n)=x_n$ the system was in at time $n$, and not on its states $X(k)=x_k$ at previous times $k=1,2,\dots,n-1$. In other words:
\begin{equation}\label{eq:markov-condition-precise}
  q(\cdot\vert\ntuple{x}{n})
  =q_n(\cdot\vert x_n),
  \quad\vtuple{x}{n}\in\states^n,\,n=1,\dots,N-1,
\end{equation}
where $q_n(\cdot\vert x_n)$ is some probability mass function on~$\states$. The Markov chain may be non-stationary, as the transition probabilities on the right-hand side in Eq.~\eqref{eq:markov-condition-precise} are allowed to depend explicitly on the time $n$. Figure~\ref{fig:markov-probability} gives an example of a  probability tree for a Markov chain with $\states=\{a,b\}$ and $N=3$.
\par
With the local probability mass functions $m_1$ and $q_n(\cdot\vert x_n)$ we associate the linear real-valued \newconcept{expectation functionals} $\ex_1$ and $\ex_n(\cdot\vert x_n)$, given, for all real-valued maps $h$ on $\states$, by
\begin{equation}
  \ex_1(h)
  \eqdef\smashoperator{\sum_{x_1\in\states}}h(x_1)m_1(x_1)
  \quad\text{ and }\quad
  \ex_n(h\vert x_n)
  \eqdef\smashoperator{\sum_{x_{n+1}\in\states}}
  h(x_{n+1})q_n(x_{n+1}\vert x_n)
\end{equation}
Throughout, we will formulate our results using expectations, rather than probabilities.\footnote{Arguments for the `expectation approach' to probability theory were given by \citet{whittle2000}. This approach is also central in the work of \citet{finetti19745}. For classical, precise probabilities, whether we use the language of probability measures, or that of expectation operators, seems to be a matter of personal preference, as the two approaches are formally equivalent. But for the imprecise-probability models we introduce in Section~\ref{sec:towards}, it was argued by \citet{walley1991} that the (lower and upper) expectation language is mathematically superior and more expressive.} Our reasons for doing so are not merely aesthetic, or a matter of personal preference; they will become clear as we go along.
\par
In any probability tree, probabilities and expectations can be calculated very efficiently using backwards recursion.\footnote{See Chapter~3 of \citeauthor{shafer1996a}'s book~\citep{shafer1996a} on causal reasoning in probability trees, which contains a number of propositions about calculating probabilities and   expectations in probability trees. That such backwards recursion is possible, was arguably   discovered by Christiaan Huygens in the middle of the 17-th century. \citetopt[Appendix~A]{shafer1996a} discusses \citeauthor{huygens16567}'s treatment~\citep[Appendix~VI]{huygens16567} of a special case of the so-called \newconcept{Problem of Points}, where Huygens draws what is probably the first recorded probability tree, and solves the problem by backwards calculation of expectations in the tree.} Suppose that in situation~$\vtuple{x}{n}$, we want to calculate the conditional expectation $\ex(f\vert\vtuple{x}{n})$ of some real-valued map~$f$ on~$\states^N$ that may depend on the values of the states $X(1)$, \dots, $X(N)$.  Let us indicate briefly how this is done, also taking into account the simplifications due to the Markov Condition~\eqref{eq:markov-condition-precise}.
\par
For these simplifications, a prominent part will be played by the so-called \newconcept{transition operators}\footnote{The operators $\trans_n$ are also called the \newconcept{generators} of the Markov process; see \citet{whittle2000}.} $\trans_n$ and $\ttrans_n$. Consider the linear space $\allgambles(\states)$ of all real-valued maps on $\states$. Then the linear operator (transformation) $\trans_n\colon\allgambles(\states)\to\allgambles(\states)$ is defined by
\begin{equation}\label{eq:trans-linear}
  \trans_nh(x_n)
  \eqdef\ex_n(h\vert x_n)
  =\smashoperator{\sum_{x_{n+1}\in\states}}h(x_{n+1})q_n(x_{n+1}\vert x_n)
\end{equation}
for all real-valued maps $h$ on $\states$. In other words, $\trans_nh$ is the real-valued map on $\states$ whose value $\trans_nh(x_n)$ in ${x_n\in\states}$ is the conditional expectation of the random variable $h(X(n+1))$, given that the system is in state~$x_n$ at time~$n$. More generally, we also consider the linear maps~$\ttrans_n$ from $\allgambles(\states^{n+1})$ to $\allgambles(\states^n)$, defined by
\begin{equation}\label{eq:trans-linear-general}
\begin{aligned}
  \ttrans_nf\ftuple{x}{n}
  \eqdef{}&{}\trans_n(f(\ntuple{x}{n},\cdot))(x_n)\\
  ={}&{}\ex_n(f(\ntuple{x}{n},\cdot)\vert x_n)
  =\smashoperator{\sum_{x_{n+1}\in\states}}f(\ntuple{x}{n},x_{n+1})
  q_{n}(x_{n+1}\vert x_n)
\end{aligned}
\end{equation}
for all $\vtuple{x}{n}\in\states^n$ and all real-valued maps $f$ on $\states^{n+1}$.\footnote{The $\ttrans^n$ can be seen as projection operators, since (with some abuse of notation) $\ttrans_n\circ\ttrans_n=\ttrans_n$.}
\par
\begin{figure}[ht]
  \centering
  \begin{tikzpicture}[x={(1.2em,0ex)},y={(0em,1ex)},xscale=.99]
    \newcommand{\ytick}{node[rectangle,inner sep=0pt,fill,minimum width=.25em,minimum height=.5pt] {}}
  \begin{scope}[yshift=0ex]\small
    \draw[yellow!50!red!33.3333!blue!25!green!50!black,fill=yellow!50!red!33.3333!blue!25!green!50] (.5,2.71875) rectangle ++(29,-2.71875) node[yshift=1.1ex,above,pos=.5,text=black] {$\ex(f)=\ex_1(\ttrans_1\ttrans_2f)$};
    \draw (0,0) -- (30,0);
    \foreach \x in {0,30} {
      \draw[->] (\x,0) -- (\x,5);
      \foreach \y in {1,2,3,4} \draw (\x,\y) \ytick;
    }
  \end{scope}
  \begin{scope}[yshift=-8ex]\small
    \draw[yellow!50!red!50!black,fill=yellow!50!red!50] (.5,3) rectangle ++(14,-3) node[yshift=1.2ex,above,pos=.5,text=black] {\small$\ex(f\vert a)=\ttrans_1\ttrans_2f(a)$};
    \draw[blue!50!green!50!black,fill=blue!50!green!50] (15.5,1.875) rectangle ++(14,-1.875) node[yshift=1ex,above,pos=.5,text=black] {$\ex(f\vert b)=\ttrans_1\ttrans_2f(b)$};
    \draw (0,0) -- (30,0);
    \foreach \x in {0,30} {
      \draw[->] (\x,0) -- (\x,5);
      \foreach \y in {1,2,3,4} \draw (\x,\y) \ytick;
    }
  \end{scope}
  \begin{scope}[yshift=-16ex]\small
    \draw[yellow!50!black,fill=yellow!50] (.5,3.5) rectangle ++(6.5,-3.5) node[yshift=1.3ex,above,pos=.5,text=black] {$\ex(f\vert a,a)=\ttrans_2f(a,a)$};
    \draw[red!50!black,fill=red!50] (8,2.5) rectangle ++(6.5,-2.5) node[yshift=1ex,above,pos=.5,text=black] {$\ex(f\vert a,b)=\ttrans_2f(a,b)$};
    \draw[blue!50!black,fill=blue!50] (15.5,2) rectangle ++(6.5,-2) node[yshift=1ex,above,pos=.5,text=black] {$\ex(f\vert b,a)=\ttrans_2f(b,a)$};
    \draw[green!50!black,fill=green!50] (23,1.75) rectangle ++(6.5,-1.75) node[yshift=1ex,above,pos=.5,text=black] {$\ex(f\vert b,b)=\ttrans_2f(b,b)$};
    \draw (0,0) -- (30,0);
    \foreach \x in {0,30} {
      \draw[->] (\x,0) -- (\x,5);
      \foreach \y in {1,2,3,4} \draw (\x,\y) \ytick;
    }
  \end{scope}
  \begin{scope}[yshift=-24ex]\small
    \draw[yellow!75!black,fill=yellow!75] (.5,4) rectangle ++(2.75,-4) node[yshift=1.6ex,above,pos=.5,text=black] {$f(a,a,a)$};
    \draw[yellow!25!black!50,fill=yellow!25] (4.25,3) rectangle ++(2.75,-3) node[yshift=1ex,above,pos=.5,text=black] {$f(a,a,b)$};
    \draw[red!75!black,fill=red!75] (8,2) rectangle ++(2.75,-2) node[yshift=1ex,above,pos=.5,text=black] {$f(a,b,a)$};
    \draw[red!25!black!50,fill=red!25] (11.75,3) rectangle ++(2.75,-3) node[yshift=1ex,above,pos=.5,text=black] {$f(a,b,b)$};
    \draw[blue!75!black,fill=blue!75] (15.5,2.5) rectangle ++(2.75,-2.5) node[yshift=1ex,above,pos=.5,text=black] {$f(b,a,a)$};
    \draw[blue!25!black!50,fill=blue!25] (19.25,1.5) rectangle ++(2.75,-1.5) node[yshift=1ex,above,pos=.5,text=black] {$f(b,a,b)$};
    \draw[green!75!black,fill=green!75] (23,.5) rectangle ++(2.75,-.5) node[yshift=1ex,above,pos=.5,text=black] {$f(b,b,a)$};
    \draw[green!25!black!50,fill=green!25] (26.75,3) rectangle ++(2.75,-3) node[yshift=1ex,above,pos=.5,text=black] {$f(b,b,b)$};
    \draw (0,0) -- (30,0);
    \foreach \x in {0,30} {
      \draw[->] (\x,0) -- (\x,5);
      \foreach \y in {1,2,3,4} \draw (\x,\y) \ytick;
    }
  \end{scope}
\end{tikzpicture}
\caption{Backwards calculation of the conditional and joint expectations of a real-valued map~$f$ on~$\states^3$, for a stationary Markov chain with state set $\states=\{a,b\}$, and a uniform probability mass function attached to each non-terminal situation.}
\label{fig:communicating-vessels}
\end{figure}
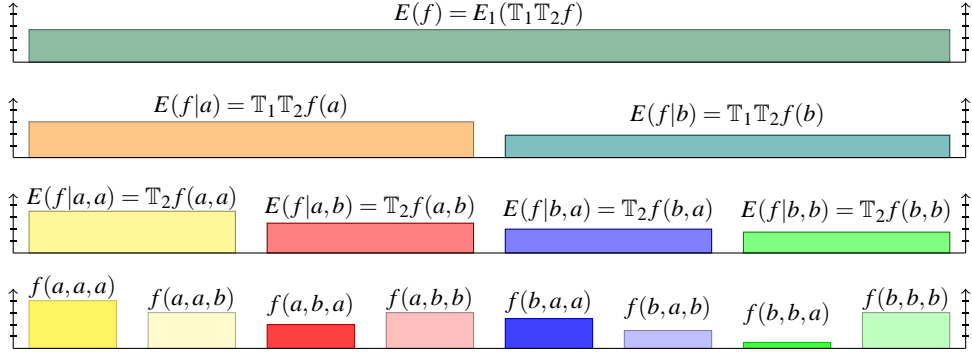

We begin our illustration of backwards recursion by calculating $\ex(f\vert\ntuple{x}{n})$ for the case $n=N-1$. Here
\begin{align}
  \ex(f\vert\ntuple{x}{N-1})
  &=\ex(f(\ntuple{x}{N-1},\cdot)\vert\ntuple{x}{N-1})\notag\\
  &=\smashoperator{\sum_{x_N\in\states}}f(\ntuple{x}{N-1},x_N)q(x_N\vert\ntuple{x}{N-1})\notag\\
  &=\smashoperator{\sum_{x_N\in\states}}f(\ntuple{x}{N-1},x_N)q_{N-1}(x_N\vert x_{N-1})
  =\ttrans_{N-1}f\ftuple{x}{N-1},
\end{align}
where the third inequality follows from the Markov Condition~\eqref{eq:markov-condition-precise}, and the fourth from Eq.~\eqref{eq:trans-linear-general}. Using similar arguments for $n=N-2$, we derive from the Law of Iterated Expectations\footnote{Also known as the Rule of Total Expectation, or the Rule of Total Probability, or the Conglomerative Property; see, e.g., \citetopt[Section~5.3]{whittle2000} or \citet{finetti19745}.} that
\begin{equation}
  \ex(f\vert\ntuple{x}{N-2})
  =\ex(\ex(f(\ntuple{x}{N-2},\cdot,\cdot)\vert\ntuple{x}{N-2},\cdot)
  \vert\ntuple{x}{N-2})
  =\ttrans_{N-2}\ttrans_{N-1}f\ftuple{x}{N-2}.
\end{equation}
Repeating this argument leads to the backwards recursion formulae
\begin{equation}\label{eq:backpropagation-precise-1}
  \ex(f\vert\ntuple{x}{n})
  =\ttrans_n\ttrans_{n+1}\dots\ttrans_{N-1}f\ftuple{x}{n}
\end{equation}
for $n=1,\dots,N-1$, while for $n=0$, we get
\begin{equation}\label{eq:backpropagation-precise-2}
  \ex(f)\eqdef\ex(f\vert\init)=\ex_1(\ttrans_1\ttrans_2\dots\ttrans_{N-1}f).
\end{equation}
In these formulae, $f$ is any real-valued map on $\states^N$.
In Figure~\ref{fig:communicating-vessels}, we give a graphical representation of calculations using the backwards recursion formulae~\eqref{eq:backpropagation-precise-1} and~\eqref{eq:backpropagation-precise-2}, for a two-state stationary Markov chain.
\par
\par
For instance, if we let $f$ be the indicator functions $\ind{\{\vtuple{x}{N}\}}$ of the singletons $\{\vtuple{x}{N}\}$, Formulae \eqref{eq:backpropagation-precise-1} and~\eqref{eq:backpropagation-precise-2} allow us to calculate the joint probability mass function $p$ defined by $p\ftuple{x}{N}=\ex(\ind{\{\vtuple{x}{N}\}})$ for all the variables $X(1)$, \dots, $X(N)$.
We can also use them to find the conditional mass functions $p_n(\cdot\vert x_n)$ and $p(\cdot\vert\ntuple{x}{n})$ defined by $p_n(\ntuple[n+1]{x}{N}\vert x_n)=p(\ntuple[n+1]{x}{N}\vert\ntuple{x}{n})=\ex(\ind{\{\vtuple{x}{N}\}}\vert\ntuple{x}{n})$.

\subsection{The Perron--Frobenius Theorem for classical Markov chains}
We are especially interested in the case of a \newconcept{stationary} Markov chain, and in the (marginal) expectation $\ex_n(h)$ of a real-valued map $h$ (on $\states$) that depends only on the state $X(n)$ at time~$n$. Here, Eq.~\eqref{eq:backpropagation-precise-2} becomes
\begin{equation}\label{eq:backpropagation-precise-3}
   \ex_n(h)\eqdef\ex_1(\trans^{n-1}h),
\end{equation}
where $\trans\eqdef\trans_1=\trans_2=\dots=\trans_{N-1}$, and where we denote by~$\trans^k$ the $k$-fold composition of $\trans$ with itself; in particular, $\trans^0$ is the identity operator $\id$ on $\allgambles(\states)$. If we let $h=\ind{\{x_n\}}$, this allows us to find the probability mass function $m_n(x_n)=\ex_n(\ind{\{x_n\}})$, $x_n\in\states$ for the state $X(n)$. 
\par
By the way, the linear transition operator $\trans$ is very closely related to the so-called \newconcept{Markov}, or \newconcept{transition}, \newconcept{matrix} $\transmat$ of the stationary Markov chain, whose elements for all $(x,y)\in\states^2$ are defined by
\begin{equation}
  T_{xy}
  \eqdef 
  q(y\vert x)=\trans\ind{\{y\}}(x).
\end{equation}
Any such transition matrix satisfies the conditions $T_{xy}\geq0$ and $\sum_{z\in\states}T_{xz}=1$. We will henceforth call \newconcept{transition matrix} any matrix satisfying these properties.\footnote{In the literature we also find the term \newconcept{stochastic matrix}, see \citet{hartfiel1998}, for instance.}
The probability counterpart of the expectation formula~\eqref{eq:backpropagation-precise-3} can then be written in matrix form as:
\begin{equation}\label{eq:back-propagation-precise-4}
  \dismat_n=\dismat_1T^{n-1},
\end{equation}
where, here and further on, we also use the notation $\dismat_n$ for the row vector whose components are the probabilities $\dismat_n(x_n)$, $x_n\in\states$. 
\par
Under some restrictions on the transition operator $\trans$, the classical Perron--Frobenius Theorem then tells us that, as $n$ (as well as the time horizon $N$) recedes to infinity, this probability mass function $m_n$ converges to some limit, independently of the initial probability mass function $m_1$; see \citetopt[Theorem~4.1.6]{kemeny1976} and \citetopt[Chapter~6]{luenberger1979}. In terms of expectation functionals and transition operators:

\begin{theorem}[Classical Perron--Frobenius Theorem, Expectation Form]\label{theo:perron-frobenius-classical}
Consider a stationary Markov chain with finite state set $\states$ and transition operator $\trans$.
Suppose that\/ $\trans$ is regular, meaning that there is some $k>0$ such that ${\min\trans^k\ind{\{x\}}>0}$ for all~$x$ in~$\states$.\footnote{This means that there is a $k>0$ such that all elements of the $k$-th power $\transmat^k$ of the transition matrix $\transmat$ are (strictly) positive. Matrices with this property are sometimes called \newconcept{regular} as well, but this same name is also used for other matrix properties. Another name for this property is `\newconcept{primitive}' \cite{hartfiel1998}.}
Then for every initial expectation operator $\ex_1$, the expectation operator $\ex_n=\ex_1\circ\trans^{n-1}$ for the state at time $n$ converges point-wise to the same limit expectation operator $\ex_\infty$:
\begin{equation}
  \smashoperator{\lim_{n\to\infty}}\ex_n(h)
  =\smashoperator{\lim_{n\to\infty}}\ex_1(\trans^{n-1}h) \defeq \ex_\infty(h)
  \quad\text{ for all $h\in\allgambles(\states)$}.
\end{equation}
Moreover, the limit expectation~$\ex_\infty$ is the only $\trans$-invariant expectation on $\allgambles(\states)$, in the sense that $\ex_\infty=\ex_\infty\circ\trans$.
\end{theorem}

\section{Towards imprecise Markov chains}\label{sec:towards}
The treatment above rests on the assumption that the initial probabilities and the transition probabilities are precisely known. If such is not the case, then it seems necessary to perform some kind of sensitivity analysis, in order to find out to what extent any conclusions we might reach using such a treatment, depend on the actual values of these probabilities.
\par
A very general way of performing a sensitivity analysis for probabilities involves calculations with closed convex sets of probability mass functions, also called \newconcept{credal sets}, rather than with single probability measures.
Let~$\simplex_\states$ denote the set of all probability mass functions on~$\states$, an~$(\abs{\states}-1)$-dimensional unit simplex in the $\abs{\states}$-dimensional linear space $\reals^\states$, then $\set{m\in\simplex_\states}{(\forall x\in\states) m(x)\leq\frac{1}{2} }$ is a cre\-dal set, but $\set{m\in\simplex_\states}{(\exists x\in\states) m(x)\geq\frac{1}{2} }$ is not.
\par
There is a growing body of literature on this interesting and fairly new area of \newconcept{imprecise probabilities}, starting with the publication of \citeauthor{walley1991}'s \citep{walley1991} seminal work.
We refer to the literature \Citep{walley1991,walley1996,weichselberger2001,cooman2005c} for more details and discussion.
\par
Let us recall a number of results for credal sets, important for the developments in this paper. Proofs can be found in \citeauthor{walley1991}'s book \citep[Chapters~2 and~3]{walley1991}.
Specifying a closed convex set~$\mass$ of probability mass functions~$p$ on a finite set~$\pties$ is equivalent to specifying its \newconcept{lower} and \newconcept{upper expectation} (functionals) $\lex_\mass\colon\allgambles(\pties)\to\reals$ and $\uex_\mass\colon\allgambles(\pties)\to\reals$, defined for all $g\in\allgambles(\pties)$ by
\begin{equation}\label{eq:mass-to-luex}
  \lex_\mass(g)\eqdef\min\set{\ex_p(g)}{p\in\mass}
  \quad\text{ and }\quad
  \uex_\mass(g)\eqdef\max\set{\ex_p(g)}{p\in\mass},
\end{equation}
where $\ex_p(g)=\sum_{\pty\in\pties}g(y)p(y)$ is the expectation of $g$ associated with the probability mass function $p$. In a sensitivity analysis, such functionals are quite useful, because they give tight lower and upper bounds on the expectation of any real-valued map. Since the functionals $\lex_\mass$ and $\uex_\mass$ are \newconcept{conjugate} in the sense that $\lex_\mass(g)=-\uex_\mass(-g)$ for all real-valued maps $g$ on $\pties$, one is completely determined if the other is known. Below, we concentrate on upper expectations. Any upper expectation $\uex=\uex_\mass$ associated with some credal set $\mass$ satisfies the following properties \citep[see, e.g.][Section~2.6.1]{walley1991}:
{\renewcommand\theenumi{$\uex$\ensuremath{\arabic{enumi}}}
\begin{enumerate}
\item $\min g\leq\uex(g)\leq\max g$ for all $g$ in $\allgambles(\pties)$ (boundedness);\label{eq:uex1}
\item $\uex(g_1+g_2)\leq\uex(g_1)+\uex(g_2)$ for all $g_1$ and $g_2$  in $\allgambles(\pties)$ (subadditivity);\label{eq:uex2}
\item $\uex(\lambda g)=\lambda\uex(g)$ for all real $\lambda\geq0$ and all  $g$  in $\allgambles(\pties)$ (non-negative homogeneity);\label{eq:uex3}
\item $\uex(g+\mu\cg)=\uex(g)+\mu$ for all real $\mu$ and all $g$ in $\allgambles(\pties)$ (constant additivity);\label{eq:uex4}
\item if $g_1\leq g_2$ then $\uex(g_1)\leq\uex(g_2)$ for all $g_1$ and $g_2$  in $\allgambles(\pties)$ (monotonicity);\label{eq:uex5}
\item if $g_n\to g$ point-wise then $\uex(g_n)\to\uex(g)$ for all sequences $g_n$ in $\allgambles(\pties)$ (continuity);\label{eq:uex6}
\item $\uex(g)\geq-\uex(-g)=\lex(g)$ for all $g$ in $\allgambles(\pties)$ (upper--lower consistency).\label{eq:uex7}
\end{enumerate}}
\noindent
Conversely, for any real functional~$\uex$ that is defined on $\allgambles(\pties)$ and that satisfies the conditions~\eqref{eq:uex1}--\eqref{eq:uex3}, there is a unique credal set $\mass\subseteq\simplex_\states$ such that $\uex$ coincides with the upper expectation $\uex_\mass$, namely $\mass=\set{p\in\simplex_\pties}{(\forall f\in\allgambles(\pties))\ex_p(f)\leq\uex(f)}$. Such an $\uex$ therefore automatically also satisfies conditions~\eqref{eq:uex4}--\eqref{eq:uex7}. It therefore make sense to {\it call \emph{upper expectation} any real functional $\uex$ on $\allgambles(\pties)$ that satisfies properties \eqref{eq:uex1}--\eqref{eq:uex3}.}
\par
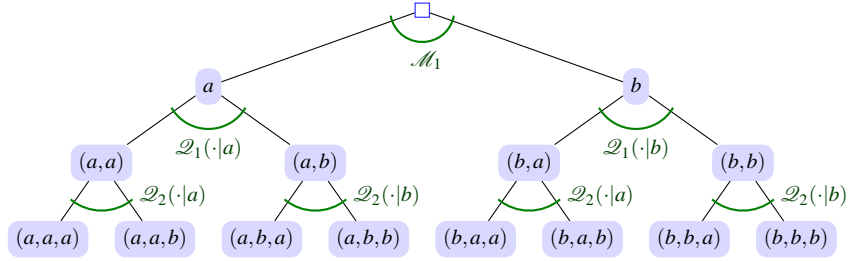
\begin{figure}[ht]
  \centering\footnotesize
  \begin{tikzpicture}
    \tikzstyle{level 1}=[sibling distance=20em]
    \tikzstyle{level 2}=[sibling distance=10em]
    \tikzstyle{level 3}=[sibling distance=5em]
    \node[root] (root) {} [grow=down,level distance=8ex]
    child {node[nonterminal] (a) {$a\vphantom{)}$}
      child {node[nonterminal] (aa) {$(a,a)$}
        child {node[nonterminal] (aaa) {$(a,a,a)$}}
        child {node[nonterminal] (aab) {$(a,a,b)$}}
      }
      child {node[nonterminal] (ab) {$(a,b)$}
        child {node[nonterminal] (aba) {$(a,b,a)$}}
        child {node[nonterminal] (abb) {$(a,b,b)$}}
      }
    }
    child {node[nonterminal] (b) {$b\vphantom{)}$}
      child {node[nonterminal] (ba) {$(b,a)$}
        child {node[nonterminal] (baa) {$(b,a,a)$}}
        child {node[nonterminal] (bab) {$(b,a,b)$}}
      }
      child {node[nonterminal] (bb) {$(b,b)$}
        child {node[nonterminal] (bba) {$(b,b,a)$}}
        child {node[nonterminal] (bbb) {$(b,b,b)$}}
      }
    };
    \draw[local,thick] (root) +(190:1.5em) arc (190:350:1.5em);
    \draw[local,thick] (b) +(210:2em) arc (210:330:2em);
    \draw[local,thick] (a) +(210:2em) arc (210:330:2em);
    \draw[local,thick] (bb) +(230:2.25em) arc (230:310:2.25em);
    \draw[local,thick] (ba) +(230:2.25em) arc (230:310:2.25em);
    \draw[local,thick] (ab) +(230:2.25em) arc (230:310:2.25em);
    \draw[local,thick] (aa) +(230:2.25em) arc (230:310:2.25em);
    \path (root) +(275:2.35em) node[local] {$\margmass_1$};
    \path (a) +(270:2.95em) node[local] {$\condmass_1(\cdot\vert a)$};
    \path (b) +(270:2.95em) node[local] {$\condmass_1(\cdot\vert b)$};
    \path (aa) +(300:2.8em) node[local,above right]  
    {$\condmass_2(\cdot\vert a)$};
    \path (bb) +(300:2.8em) node[local,above right]  
    {$\condmass_2(\cdot\vert b)$};
    \path (ba) +(300:2.8em) node[local,above right]  
    {$\condmass_2(\cdot\vert a)$};
    \path (ab) +(300:2.8em) node[local,above right]  
    {$\condmass_2(\cdot\vert b)$};
  \end{tikzpicture}
  \caption{The tree for the time evolution of an imprecise Markov chain that can be in two states, $a$ and $b$, and can change state at each time instant $n=1,2$.}
  \label{fig:markov-imprecise}
\end{figure}

What is the upshot of all this for the Markov chain problem we are considering here? First of all, in the initial situation~$\init$, corresponding to time $n=0$, rather than a single initial probability mass function~$m_1$, we now have a local credal set $\margmass_1$ of candidate mass functions~$m_1$ for the state $X(1)$ that the system will be in at time $k=1$.
We denote by $\uex_1$ the upper expectation associated with $\margmass_1$:
\begin{equation}
  \uex_1(h)
  \eqdef\max\Bigl\{\smashoperator[r]{\sum_{x\in\states}}h(x)m_1(x)
  \colon m_1\in\margmass_1\Bigr\}
  \quad\text{ for all $h\in\allgambles(\states)$.}
\end{equation}
Also, in any situation $\vtuple{x}{n}\in\states^n$ corresponding to time $n=1,2,\dots,N-1$, instead of a single transition mass function $q_n(\cdot\vert x_n)$, we now have a local credal set $\condmass_n(\cdot\vert x_n)$ of candidate conditional mass functions $q_n(\cdot\vert x_n)$ for the state $X(n+1)$ that the system will be in at time $n+1$. We denote by $\uex_n(\cdot\vert x_n)$ the upper expectation associated with $\condmass_n(\cdot\vert x_n)$, i.e.:
\begin{equation}\label{eq:local-upper}
  \uex_n(h\vert x_n)
  \eqdef\max\Bigl\{\smashoperator[r]{\sum_{x\in\states}}
    h(x)q(x)\colon q\in\condmass_n(\cdot\vert x_n) \Bigr\}
    \quad\text{ for all $h\in\allgambles(\states)$}.
\end{equation}
We call the resulting model an \newconcept{imprecise Markov chain}. Figure~\ref{fig:markov-imprecise} gives an example of a probability tree for an imprecise Markov chain.
It is an imprecise-probability tree where the local conditional models satisfy the \emph{Markov Condition}:
\begin{equation}\label{eq:imprecise-markov-condition}
  \condmass(\cdot\vert\ntuple{x}{n})
  =\condmass(\cdot\vert x_n)
  \quad\text{ for all $\ntuple{x}{n}\in\states^n$ and $n=1,2,\dots,N-1$}.
\end{equation}
A classical, or \newconcept{precise}, Markov chain is an imprecise one with credal sets that are singletons.
\par
\par
How, then, can a sensitivity analysis be performed for such an imprecise Markov chain? We choose, in each non-terminal situation $\vtuple{x}{k}$ of the above-mentioned event tree, a local transition probability mass $q(\cdot\vert\ntuple{x}{k})$ in the set of possible candidates $\condmass_k(\cdot\vert x_k)$.\footnote{These local transition probability masses themselves depend on the situation $\vtuple{x}{k}$ they are attached to, but the sets $\condmass_k(\cdot\vert x_k)$ they are chosen from only depend on the last state~$x_k$, as the Markov Condition~\eqref{eq:imprecise-markov-condition} tells us.} For $k=0$, we get the initial situation~$\init$, where we choose some element~$m_1$ in the set of possible candidates~$\margmass_1$.
By making a choice of local model for each non-terminal situation in the event tree, we obtain what we call a \newconcept{compatible probability tree}, for which we may calculate all (conditional) expectations and probability mass functions:
\begin{align}
  \ex(f\vert\ntuple{x}{n})
  &=\smashoperator{\sum_{\vtuple[n+1]{x}{N}\in\states^{N-n}}}f(\ntuple{x}{n},\ntuple[n+1]{x}{N})
  \smashoperator{\prod_{k=n}^{N-1}}q(x_{k+1}\vert\ntuple{x}{k})
  \label{eq:probability-tree-conditional},\\
  \ex(f)
  &=\smashoperator{\sum_{\vtuple{x}{N}\in\states^{N}}}f\ftuple{x}{N}
  m_1(x_1)\smashoperator{\prod_{k=1}^{N-1}}q(x_{k+1}\vert\ntuple{x}{k}),
  \label{eq:probability-tree-joint}
\end{align}
for ${n=1,\dots,N-1}$, and for all real-valued maps~$f$ on~$\states^N$. 
As we have just come to realise, the probability trees that are compatible with an imprecise Markov chain are no longer necessarily (precise) Markov chains themselves. 
It is still possible to calculate the $\ex(f\vert\ntuple{x}{n})$ and $\ex(f)$ in Eqs.~\eqref{eq:probability-tree-conditional} and~\eqref{eq:probability-tree-joint} using backwards recursion \citep[Chapter~3]{shafer1996a}, but the formulae for doing so will be more complicated than the ones for precise Markov chains given by Eqs.~\eqref{eq:backpropagation-precise-1} and~\eqref{eq:backpropagation-precise-2}.
\par
If we repeat this for every other choice of the $m_1$ in $\margmass_1$ and the $q(\cdot\vert\ntuple{x}{k})$ in $\condmass_k(\cdot\vert x_k)$, we end up with an infinity of compatible probability trees,\footnote{Except when all the credal sets are singletons, of course.} for which the associated (conditional) expectations and probability mass functions turn out to constitute closed convex sets.
We denote their corresponding upper expectation functionals on $\allgambles(\states^N)$ by $\uex(\cdot\vert\ntuple{x}{n})$ and $\uex$. These upper expectations, and the conjugate lower expectations, are the final aim of our sensitivity analysis.
\par
The procedure we have just described is computationally very complex. When the closed convex sets $\margmass_1$ and $\condmass_k(\cdot\vert x)$ each have a finite number of extreme points (are polytopes), we can limit ourselves to working with these sets of extreme points, rather than with the infinite sets themselves. But even then, the computational complexity of this approach will generally  be exponential in the number of time steps.
\par
However, we will see in Section~\ref{sec:sensitivity-analysis} that the upper expectations $\uex$ and $\uex(\cdot\vert\ntuple{x}{n})$ associated with the closed convex sets of (conditional) probability mass functions for the compatible probability trees of an imprecise Markov chain can be calculated in the same way as the expectations $\ex$ and $\ex(\cdot\vert\ntuple{x}{n})$ in a precise one: using counterparts of the backwards recursion formulae~{\eqref{eq:backpropagation-precise-1}--\eqref{eq:backpropagation-precise-3}}. Because of this, making inferences about the mass function of the state at time $n$, i.e., finding the upper envelope $\uex_n$ of the $\ex_n$ given in Eq.~\eqref{eq:backpropagation-precise-3} \emph{now has a complexity that is linear, rather than exponential, in the number of time steps $n$.} This is our first contribution.
\par
Our second contribution in this paper is a Perron--Frobenius Theorem for a special class of so-called regularly absorbing stationary imprecise Markov chains: in Section~\ref{sec:convergence} we prove a generalisation of Theorem~\ref{theo:perron-frobenius-classical}, which tells us that under fairly weak conditions, the upper expectation operators~$\uex_n$ converge to limits that do not depend on the initial upper expectation operators~$\uex_1$. Our result also extends a number of other related convergence theorems for imprecise Markov chains in the literature \citep{hartfiel1991,hartfiel1994,hartfiel1998,skulj2007}.

\section{Sensitivity analysis of imprecise Markov chains}\label{sec:sensitivity-analysis}
We can now take our most important step: deriving the backwards recursion formulae for the conditional and joint upper expectations in an imprecise Markov chain.  
We first define \newconcept{upper transition operators} $\utrans_n$ and $\uttrans_n$.  
The operator $\utrans_n\colon\allgambles(\states)\to\allgambles(\states)$ is defined by
\begin{equation}\label{eq:trans-upper}
  \utrans_nh(x_n)\eqdef\uex_n(h\vert x_n)
\end{equation}
for all real-valued maps $h$ on $\states$, and all $x_n$ in $\states$. In other words, $\utrans_nh$ is the real-valued map on $\states$, whose value $\utrans_nh(x_n)$ in $x_n\in\states$ is the conditional upper expectation of the random variable $h(X(n+1))$, given that the system is in state $x_n$ at time $n$. More generally, we also consider the maps $\uttrans_n$ from $\allgambles(\states^{n+1})$ to $\allgambles(\states^n)$, defined by
\begin{equation}\label{eq:trans-upper-general}
  \uttrans_nf\ftuple{x}{n} \eqdef \bigl(\utrans_nf(\ntuple{x}{n},\cdot)\bigr)(x_n) = \uex_n(f(\ntuple{x}{n},\cdot)\vert x_n)
\end{equation}
for all $\vtuple{x}{n}$ in $\states^n$ and all real-valued maps $f$ on $\states^{n+1}$. Of course, we can also consider lower expectations and lower transition operators, which are related to the upper expectations and upper transition operators by conjugacy. As is the case for upper expectations, it is possible to introduce the notion of an upper transition operator directly, by basing it on a number of defining properties, rather than by referring to an underlying imprecise Markov chain. We refer to the Appendix for more details.
\par
The upper expectations $\uex(\cdot\vert\ntuple{x}{n})$ and $\uex$ on $\allgambles(\states^N)$ can be calculated very easily by backwards recursion, cfr.~\eqref{eq:backpropagation-precise-1} and~\eqref{eq:backpropagation-precise-2}.

\begin{theorem}[Concatenation Formula]\label{theo:concatenation}
  For any $\vtuple{x}{n}$ in $\states^n$, $n=1,\dots,N-1$, and for any real-valued map~$f$ on~$\states^N$:
\begin{align}
  \uex(f\vert\ntuple{x}{n})
  &=\uttrans_n\uttrans_{n+1}\dots\uttrans_{N-1}f\ftuple{x}{n}
  \label{eq:backpropagation-upper-1}\\
  \uex(f)
  &=\uex_1(\uttrans_1\uttrans_2\dots\uttrans_{N-1}f).
  \label{eq:backpropagation-upper-2}
\end{align}
\end{theorem}

Call, for any non-empty subset $I$ of $\{1\dots,N\}$, a real-valued map $f$ on $\states^N$ \newconcept{\mbox{$I$-measurable}} if $f\ftuple{x}{N}=f\ftuple{z}{N}$ for all $\vtuple{x}{N}$ and $\vtuple{z}{N}$ in $\states^N$ such that $x_k=z_k$ for all $k\in I$.
In other words, an $I$-measurable~$f$ only depends on the states $X(k)$ at times $k\in I$.
As an example, an \mbox{$\{n\}$-measurable} map~$h$ only depends on the state $X(n)$ at time~$n$, and we identify it with a map on~$\states$ (but remember that it acts on states at time~$n$).
The following proposition tells us that all conditional upper expectations satisfy a Markov Condition  (cfr.~\eqref{eq:markov-condition-precise}).

\begin{proposition}[Markov Condition]\label{prop:markov}
  Consider an imprecise Markov chain with finite state set $\states$ and time horizon $N$. 
Fix~$n\in\{1,\dots,N-1\}$. 
Let~$\vtuple{x}{n-1}$ and~$\vtuple{z}{n-1}$ be arbitrary elements of~$\states^{n-1}$, and let $x_n\in\states$. 
Let~$f$ be any $\{n,n+1,\dots,N\}$-mea\-sur\-able real-valued map on $\states^N$. 
Then  $\uex(f\vert\ntuple{x}{n-1},x_n)=\uex(f\vert\ntuple{z}{n-1},x_n)$, so we may write
$\uex(f\vert\ntuple{x}{n-1},x_n)=\uex_{\vert n}(f\vert x_n)$.
\end{proposition}
\noindent
The index `$\vert n$' is intended to make clear that we are considering an expectation conditional on the state $X(n)$ at time $n$.
\par
If we apply the joint upper expectation~$\uex$ to maps~$h$ that only depend on the state $X(n)$ at time~$n$, we get the \newconcept{marginal upper expectation} $\uex_n(h)\eqdef\uex(h)$, and $\uex_n$ is a model for the uncertainty about the state $X(n)$ at time~$n$. More generally, taking into account Proposition~\ref{prop:markov}, we use the notation $\uex_{n\vert\ell}(h\vert x_\ell)\eqdef\uex_{\vert\ell}(h\vert x_\ell)$ for the upper expectation of $h(X(n))$, conditional on $X(\ell)=x_\ell$ with $1\leq\ell<n$.  With notations established in Eq.~\eqref{eq:local-upper}, $\uex_{n+1\vert n}(h\vert x_n)=\uex_n(h\vert x_n)=\utrans_nh(x_n)$. Such expectations can be found using simpler recursion formulae than Eqs.~\eqref{eq:backpropagation-upper-1} and~\eqref{eq:backpropagation-upper-2}, as they are based on the simpler upper transition operators $\utrans_k$.

\begin{corollary}\label{cor:marginal-concatenation}
  For any real-valued map~$h$ on~$\states$, and for any $1\leq\ell<n\leq N$ and all~$x_\ell$ in~$\states$:
  \begin{equation}\label{eq:backpropagation-upper-3}
    \uex_{n\vert\ell}(h\vert x_\ell)
    =\utrans_\ell\utrans_{\ell+1}\dots\utrans_{n-1}h(x_\ell)
    \quad\text{ and }\quad
    \uex_n(h)
    =\uex_1(\utrans_1\utrans_2\dots\utrans_{n-1}h).
  \end{equation}
\end{corollary}
\noindent
This offers a reason for formulating our theory in terms of real-valued maps rather than events: suppose we want to calculate the upper probability $\uex_n(A)$ that the state $X(n)$ at time~$n$ belongs to the set~$A$.
According to Eq.~\eqref{eq:backpropagation-upper-3}, $\uex_n(A)=\uex_1(\utrans_1\dots\utrans_{n-1}\ind{A})$, and even if~$\utrans_{n-1}\ind{A}$ can still be calculated using upper probabilities only, it will generally assume values other than~$0$ and~$1$, and therefore will generally not be the indicator of some event.
Already after one step, i.e., in order to calculate $\utrans_{n-2}\utrans_{n-1}\ind{A}$, we need to leave the ambit of events, and turn to the more general real-valued maps; even if we only want to calculate upper \emph{probabilities} after~$n$ steps.

For joint upper and lower probability mass functions, however, we can remain within the ambit of events:

\begin{proposition}[Chapman--Kolmogorov Equations]\label{prop:chapman-kolmogorov}
  For an imprecise Markov chain, we have for all\/ $1\leq n<m\leq N$ and all\/ $(x_n,\ntuple[n+1]{x}{m})\in\states^{m-n+1}$ that
  \begin{equation}\label{eq:CKu}
    \uex_{\vert n}(\{\vtuple[n+1]{x}{m}\}\vert x_n)
    =\smashoperator{\prod_{k=n}^{m-1}}\utrans_k\ind{\{x_{k+1}\}}(x_k),
  \end{equation}
  and for all\/ $1\leq m\leq N$ and all $\vtuple{x}{m}\in\states^{m}$ that
  \begin{equation}\label{eq:jCKu}
    \uex(\{\vtuple{x}{m}\})
    =\uex_1(\{x_1\})\smashoperator{\prod_{k=1}^{m-1}}
    \utrans_k\ind{\{x_{k+1}\}}(x_k).
  \end{equation}
There are analogous expressions for the lower expectations.
\end{proposition}

\section{Accessibility relations}\label{sec:accessibility}
From now on, and for the rest of the paper, we mainly consider \emph{stationary imprecise Markov chains with an infinite time horizon}.
This means that for each time $n\in\naturals$, we consider the same upper transition operator $\utrans_n=\utrans$.

The classification of the states of such a stationary (im)precise Markov chain can be fruitfully started by introducing a so-called \newconcept{accessibility relation} $\access[\cdot]{\cdot}{\cdot}$: let~$x$ and~$y$ be any two states in~$\states$ and let~$n$ be a number of steps in~$\naturals_0\eqdef\naturals\cup\{0\}$, then $\access[n]{x}{y}$ expresses that~$y$ is accessible from~$x$ in~$n$ steps.
To be an accessibility relation, a generic ternary relation $\access[\cdot]{\cdot}{\cdot}$ has to satisfy the defining properties:
\begin{align}
  (\forall x,y\in\states) &\access[0]{x}{y}\asa x=y \label{eq:basic-communication-1},\\
  (\forall x,y,z\in\states) (\forall m,n\in\naturals_0) &\text{$\access[n]{x}{y}$ and $\access[m]{y}{z}$}\then\access[n+m]{x}{z} \label{eq:basic-communication-2}.\\
  (\forall x\in\states) (\forall n\in\naturals) (\exists y\in\states) &\access[n]{x}{y}. \label{eq:basic-communication-3}
\end{align}
\par
An accessibility relation is classically derived from the transition matrix of a stationary Markov chain; in Section~\ref{sec:accessibility-imprecise} we will associate such a relation with a stationary imprecise Markov chain.
But for \emph{any} (abstract) accessibility relation satisfying the conditions~\eqref{eq:basic-communication-1}--\eqref{eq:basic-communication-3}, we can draw all the following conclusions, no matter what transition matrix or operator it was derived from, or whether it comes about in any other way; \citetopt[Section~1.4]{kemeny1976} give a detailed justification.
In what follows, we use the terminology introduced by \citeauthor{kemeny1976}, but we want to remind the reader that the terms we use may also have various other meanings in different parts of the literature.

\subsection{Abstract accessibility relations}\label{sec:accessibility-abstract}
\begin{figure}[htb]
  \centering
  \begin{tikzpicture}[->,>=stealth,shorten >=2pt,shorten <=2pt,node distance=2em]
    \node[statesbackground,label={17:\small$\states$}] {
      \begin{tikzpicture}
        \node[closedbackground,label={40:\small$C_1$}] (C1) {
          \tikzstyle{commun}=[communbackground,minimum height=5ex,minimum width=2.8em]
          \begin{tikzpicture}[node distance=4ex and 0em]
            \node[commun,label={4:\small$D_3$}] (D3) {};
            \node[below left=of D3,commun,label={4:\small$D_1$}] (D1) {};
            \node[below right=of D3,commun,label={4:\small$D_2$}] (D2) {};
            \node[above left=of D3,commun,label={4:\small$D_4$}] (D4) {};
            \node[above right=of D3,commun,label={4:\small$D_5$}] (D5) {};
            \draw (D1) -- (D3);
            \draw (D2) -- (D3);
            \draw (D3) -- (D4);
            \draw (D3) -- (D5);
          \end{tikzpicture}
        };
        \node[right=of C1,closedbackground,label={40:\small$C_2$}] (C2) {
          \tikzstyle{commun}=[communbackground,minimum height=10ex,minimum width=2.8em]
          \begin{tikzpicture}[node distance=4ex and 0em]
            \node[commun,label={30:\small$D_8$}] (D8) {};
            \node[below left=of D8,commun,label={30:\small$D_6$}] (D6) {};
            \node[below right=of D8,commun,label={30:\small$D_7$}] (D7) {};
            \draw (D6) -- (D8);
            \draw (D7) -- (D8);
          \end{tikzpicture}
        };
        \node[right=of C2,closedbackground,label={59:\small$C_3$}] (C3) {
          \tikzstyle{commun}=[communbackground,minimum height=24ex,minimum width=2.8em]
          \begin{tikzpicture}
            \node[commun,label={70:\small$D_9$}] (D9) {};
          \end{tikzpicture}
        };
      \end{tikzpicture}
    };
  \end{tikzpicture}
  \caption{
    Three increasingly finer partitions of the state set $\states$ for a particular stationary (im)precise Markov chain, or more generally, for an accessibility relation $\access[\cdot]{\cdot}{\cdot}$.
    No transition between states of the classes $C_1$, $C_2$, and $C_3$ is possible, and these classes can be seen as separate (im)precise Markov chains.
    The equivalence classes $D_k$ for the communication relation are partially ordered by the relation $\access{}{}$, whose (Hasse) diagram is represented by the upward arrows.
    Maximal classes are $D_4$, $D_5$, $D_8$, and $D_9$, the other classes are transient.
    If $D_4$, $D_5$, $D_8$, and $D_9$ are aperiodic, the accessibility relation restricted to respectively $C_1$, $C_2$, and $C_3$ is respectively maximal class regular, top class regular, and regular.}
  \label{fig:communication}
\end{figure}
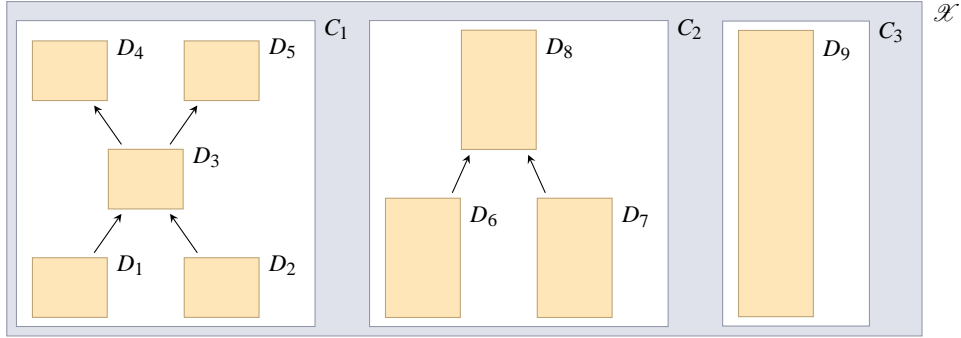

Accessibility relations give rise to many interesting concepts,  which we discuss below.
We refer to Figure~\ref{fig:communication} for a graphical representation.
\par
\par
Consider any two states $x$ and $y$ in $\states$.
Then $y$ is \newconcept{accessible from $x$}, which we denote as $\access{x}{y}$, if there is some $n\in\naturals_0$ such that $\access[n]{x}{y}$.
If $x$ and $y$ are accessible from one another, then we say that $x$ and $y$ \newconcept{communicate}, which we denote as $\commun{x}{y}$.
\par
It follows from Eqs.~\eqref{eq:basic-communication-1} and~\eqref{eq:basic-communication-2} that the binary relation $\access{}{}$ on $\states$ is a preorder, i.e., is reflexive and transitive. The binary relation $\commun{}{}$ on~$\states$ is the associated equivalence relation. This \newconcept{communication relation} $\commun{}{}$ partitions the state set~$\states$ into equivalence classes~$D$ of states that are accessible from one another, called \newconcept{communication classes}. The preorder $\access{}{}$ induces a partial order on this partition, also denoted by $\access{}{}$.
\par
Undominated or \newconcept{maximal} states with respect to the preorder $\access{}{}$ are states~$x$ such that $\access{x}{y}\then\access{y}{x}$ for any state~$y$ in~$\states$.
This means that a maximal state has access only to other maximal states in the same communication class, and to no other states.
Collections of maximal states, such as the communication classes they belong to, are also called \newconcept{maximal}.
The other states and collections of them, such as the communication classes they belong to, are called \newconcept{transient}.
If all maximal states communicate, or in other words if there is a unique maximal communication class, this class is called the \newconcept{top} class.
It is made up of those states that are accessible from any state.
\par
Consider, for any $x$ and $y$ in $\states$, the set
\begin{equation}
  \nsteps{x}{y}\eqdef\inlineset{n\in\naturals}{\access[n]{x}{y}},
\end{equation}
i.e., those numbers of steps after which $y$ is accessible from $x$.  We call the \newconcept{period} $\period{x}$ of a state $x$ the greatest common divisor of the set $\nsteps{x}{x}$, i.e., $\period{x}\eqdef\gcd\inlineset{n\in\naturals}{\access[n]{x}{x}}$. Because, by Eq.~\eqref{eq:basic-communication-2}, $\nsteps{x}{x}$ is closed under addition, we can rely on a basic number-theoretic result (see, e.g., \citetopt[Theorem~1.4.1]{kemeny1976}) which tells us that $\nsteps{x}{x}$ is, up to perhaps a finite number of initial elements, equal to the set of all multiples of $\period{x}$. 
\par
Now consider an equivalence class~$D$ of communicating states, and any two states~$x$ and~$y$ in that class. Then it is not difficult to show that they have the same period: $\commun{x}{y}\then\period{x}=\period{y}$. We denote by $\period{D}$ the common period of all elements of the equivalence class $D$. 

\begin{proposition}\label{prop:class-cycle}
  Consider arbitrary $x$ and $y$ in some maximal communication class $D$.  Then there is some $0\leq\steps{x}{y}<\period{D}$ such that $n\in\nsteps{x}{y}\then n\equiv\steps{x}{y}\pmod{\period{D}}$, i.e., $n$ and $\steps{x}{y}$ are equal up to some multiple of $\period{D}$. Moreover, 
  \begin{equation}\label{eq:nxy}
    (\exists n\in\naturals)
    (\forall k\geq n)\,
    \steps{x}{y}+k\period{D}\in\nsteps{x}{y}.
  \end{equation}
\end{proposition}
\noindent
For any~$x$, $y$~and~$z$ in this equivalence class~$D$, $\steps{x}{y}+\steps{y}{z}\equiv\steps{x}{z}\pmod{\period{D}}$, and therefore $\steps{y}{z}=0$ if and only if $\steps{x}{y}=\steps{x}{z}$. This implies that `$\steps{y}{z}=0$' determines an equivalence relation on this equivalence class~$D$, which further partitions it into~$\period{D}$ subsets, called \newconcept{cyclic classes}. In such a cyclic class, all states~$y$ give the same value to~$\steps{x}{y}$, for any given~$x$ in~$D$. Within $D$, the system moves from cyclic class to cyclic class, in a definite ordered cycle of length~$\period{D}$. If~$D$ is transient, then in some cyclic classes it is possible that, rather than moving to the next cyclic class, the system moves to (a state in) another equivalence class~$D'$ for the communication relation that is a successor to~$D$ for the partial order $\access{}{}$.
\par
If $\period{D}=1$, or in other words if $\steps{x}{y}=0$ for all $x,y\in D$, then there is only one cyclic class in $D$, and we call the communication class $D$, and all its states, \newconcept{aperiodic}. If $D$ is moreover maximal, then $D$ is called \newconcept{regular}. The following general characterisation of regularity is easily derived from Proposition~\ref{prop:class-cycle}; see also  \citeauthor{kemeny1976}'s arguments \citep[Chapters~1 and~4]{kemeny1976}.

\begin{proposition}\label{prop:regularity-characterisation}
  A communication class $D\subseteq\states$ is regular under the accessibility relation $\access[\cdot]{\cdot}{\cdot}$ if and only if
  \begin{equation}\label{eq:regular}
    (\exists n\in\naturals)
    (\forall k\geq n)
    (\forall x,y\in D)
    \access[k]{x}{y}.
  \end{equation}
\end{proposition}

An interesting special case obtains when there is only one equivalence class for the communication relation (namely~$\states$), so $\states$ is maximal, and there is only one cyclic class (namely~$\states$), meaning that all states are aperiodic. In that case, the accessibility relation $\access[\cdot]{\cdot}{\cdot}$ is called \newconcept{regular} as well.
\noindent 
If all maximal communication classes are regular (aperiodic), the accessibility relation is called \newconcept{maximal class regular}. If there is only one maximal communication class, and if this top class is moreover regular (aperiodic), then the accessibility relation is called \newconcept{top class regular}. Top class regularity has the following simple alternative characterisation.

\begin{proposition}\label{prop:topclassregular}
  An accessibility relation $\access[\cdot]{\cdot}{\cdot}$ is top class regular if and only if the corresponding set $\maxregstates_{\access[]{}{}}$ of so-called \newconcept{maximal regular states} is non-empty:
  \begin{equation}\label{eq:topclassregular}
    \maxregstates_{\access[]{}{}}
    =\set{x\in\states}
    {(\exists n\in\naturals)(\forall k\geq n)(\forall y\in\states)\access[k]{y}{x}}
    \neq\emptyset;
  \end{equation}
  and in that case this set $\maxregstates_{\access[]{}{}}$ is the top communication class.
\end{proposition}

\subsection{Accessibility relations for imprecise Markov chains}\label{sec:accessibility-imprecise}
Because we now only consider stationary imprecise Markov chains, this means that for each time $n\in\naturals$, we consider the same transition models $\condmass_n(\cdot\vert x)=\condmass(\cdot\vert x)$, $x\in\states$, or equivalently, for the upper transition operators: $\utrans_n=\utrans$ and $\uttrans_n=\uttrans$.
\par
Let us denote by $\smash[b]{\utp[n]{x}{y}}$ the upper probability of going in $n$ steps from state $x$ to state $y$.
For $n=0$, $\utp[0]{x}{y}=\ind{\{y\}}(x)$, and for $n\geq1$, $\utp[n]{x}{y}=\uex_{k+n\vert k}(\{y\}\vert x)$, where---because of stationarity---the right-hand sides does not depend on~$k\in\naturals$.
By Corollary~\ref{cor:marginal-concatenation}, we find that $\utp[n]{x}{y}=\utrans^n\ind{\{y\}}(x)$ for all $n\in\naturals_0$.
The following two propositions allow us to associate an accessibility relation with the upper transition operator. 
They are immediate generalisations of similar results involving (precise) probabilities in (precise) Markov chains:

\begin{proposition}\label{prop:basic-inequality}
  For all $x$, $y$ and $z$ in $\states$, and for all $n$ and $m$ in $\naturals_0$,
  \begin{equation}\label{eq:basic-inequality}
    \utp[n+m]{x}{y}\geq\utp[n]{x}{z}\utp[m]{z}{y}.
  \end{equation}
\end{proposition}

\begin{proposition}\label{prop:always-arrive-in-a-state}
  For all $x$ in $\states$, and for all $n$ in $\naturals_0$, there is some $y$ in $\states$ such that $\utp[n]{x}{y}>0$.
\end{proposition}
\noindent
Because of these results, which ensure that Eqs.~\eqref{eq:basic-communication-2} and~\eqref{eq:basic-communication-3} are satisfied [Eqs.~\eqref{eq:basic-communication-1} is trivially satisfied because $\utp[0]{x}{y}=\ind{\{y\}}(x)$], we can define an accessibility relation $\uaccess[\cdot]{\cdot}{\cdot}$ using the $\utp[n]{x}{y}$: for any $x$ and $y$ in $\states$ and any $n\in\naturals_0$:
\begin{equation}\label{eq:uaccessibility}
  \uaccess[n]{x}{y}\asa\utp[n]{x}{y}>0\asa\utrans^n\ind{\{y\}}(x)>0.
\end{equation}
Clearly, $\uaccess[n]{x}{y}$ if there is \emph{some} compatible probability tree in which it is possible (meaning that there is a non-zero probability) to go from state $x$ to $y$ in $n$ time steps. In other words, $\uaccess[n]{x}{y}$ if it is not considered impossible in the context of our imprecise-probability model to go from $x$ to $y$ in $n$ steps: we then say that $y$ is \newconcept{accessible} from $x$ in $n$ steps; and if $\uaccess[]{x}{y}$ then $y$ is \newconcept{accessible} from $x$.
\par
The following notion will be essential for the convergence result we present in the next section. It involves both lower and upper transition probabilities.

\begin{definition}[Regularly absorbing]\label{def:regabs}
  A stationary imprecise Markov chain is called \newconcept{regularly absorbing} if it is top class regular (under~$\uaccess{}{}$), meaning that
  \begin{equation}
    \maxregstates_{\uaccess{}{}}
    \eqdef\set{x\in\states}
    {(\exists n\in\naturals)(\forall k\geq n)(\forall y\in\states)
    \utrans^k\ind{\{x\}}(y)>0}
    \neq\emptyset,
  \end{equation}
  and if moreover for all~$y$ in~$\states\setminus\maxregstates_{\uaccess[]{}{}}$ there is some  $n\in\naturals$ such that\/ $\ltrans^n\ind{\maxregstates_{\uaccess{}{}}}(y)>0$.
\end{definition}
\noindent
In particular, an imprecise Markov chain that is regular (under~$\uaccess{}{}$, meaning that the accessibility relation $\uaccess{}{}$ is regular) is also regularly absorbing (under $\uaccess{}{}$) in a trivial way.

\section{Convergence for stationary imprecise Markov chains}\label{sec:convergence}
We call an upper expectation $\uex$ on $\allgambles(\states)$ \newconcept{$\utrans$-invariant} whenever $\uex\circ\utrans=\uex$, so whenever $\uex(\utrans h)=\uex(h)$ for all $h\in\allgambles(\states)$.

\begin{theorem}[Perron--Frobenius Theorem, Upper Expectation Form]\label{theo:convergence}
  Consider a stationary imprecise Markov chain with finite state set $\states$ that is regularly absorbing. Then for every initial upper expectation $\uex_1$, the upper expectation $\uex_n=\uex_1\circ\utrans^{n-1}$ for the state at time $n$ converges point-wise to the same upper expectation $\uex_\infty$:
\begin{equation}
  \smashoperator{\lim_{n\to\infty}}\uex_n(h)
  =\smashoperator{\lim_{n\to\infty}}
  \uex_1(\utrans^{n-1}h)\defeq\uex_\infty(h)
  \text{ for all $h$ in $\allgambles(\states)$.}
\end{equation}
Moreover, the limit upper expectation $\uex_\infty$ is the only $\utrans$-invariant upper expectation on $\allgambles(\states)$.
\end{theorem}
\noindent
Let us compare this convergence result to what exists in the literature.
\par
The classical Perron--Frobenius Theorem~\ref{theo:perron-frobenius-classical} is of course a special case of our Theorem~\ref{theo:convergence}, because if (the transition operator of) a precise stationary Markov chain is regular in the sense of Theorem~\ref{theo:perron-frobenius-classical}, then it is also regular (under~$\uaccess{}{}$), and therefore regularly absorbing.
\par
Other authors have presented convergence results for stationary imprecise Markov chains, namely \citet{hartfiel1991,hartfiel1998,hartfiel1994}, and \citet{skulj2007}. They all use the following approach. 
They consider some set $\transmats$ of (one-step) transition matrices $T$, and deduce from that a corresponding set $\transmats^n$ of $n$-step transition matrices given by
\begin{equation}
  \transmats^n
  \eqdef\set{\transmat_1\transmat_2\dots\transmat_n}
  {\transmat_1,\transmat_2,\dots,\transmat_n\in\transmats}.
\end{equation}
\citeauthor{hartfiel1998} calls the sequence $\transmats^n$, $n\in\naturals$ a \newconcept{Markov set chain}. 
If we also have a set $\margmass_1$ of (marginal) mass functions $m_1$ for $X(1)$, then they take the corresponding set $\margmass_n$ of (marginal) mass functions for $X(n)$ to be 
\begin{equation}
  \margmass_n
  =\set{\dismat_1\transmat}
  {m_1\in\margmass_1\text{ and }\transmat\in\transmats^{n-1}},
\end{equation}
where, as before, we also denote by $\dismat$ the row vector corresponding to the mass function~$m$. 
If we furthermore also denote by $h$ the column vector corresponding to the values $h(x)$ of the real-valued map~$h$ in all $x\in\states$, then we find that the corresponding set $\expects_n(h)$ of expectations of $h(X(n))$ is given by
\begin{equation}
  \expects_n(h)
  =\set{m_1Th}
  {m_1\in\margmass_1\text{ and }\transmat\in\transmats^{n-1}}.
\end{equation}
Incidentally, these are also the formulae that can be obtained by considering imprecise Markov chains to be special cases of so-called credal networks under a strong independence assumption; for more details, see \citeauthor{cozman2000}'s work \citep{cozman2000,cozman2005} for instance.
\par
\citet{skulj2007} considers the set $\transmats$ of transition matrices $\transmat$ corresponding to a so-called \newconcept{interval stochastic matrix}, meaning that $\transmats$ is the set of all transition matrices such that $\ltransmat\leq\transmat\leq\utransmat$, where $\ltransmat$ and $\utransmat$ are so-called lower and upper transition matrices; see also Section~\ref{sec:lower-upper-mass} for the related model in terms of upper transition operators. \citet{hartfiel1991} considers arbitrary sets of transition matrices, but in his book \cite{hartfiel1998} he also focuses mainly on interval stochastic matrices.
\par
What is the relationship between the Markov set-chain model and the model involving upper transition operators we have studied and motivated above? 
Consider a stationary imprecise Markov chain with upper transition operator $\utrans$. 
For each state $x$, as $\utrans h(x)$ has been defined as a conditional upper expectation $\uex(h\vert x)$, there is a corresponding credal set $\condmass_{\utrans}(\cdot\vert x)$ given by
\begin{equation}\label{eq:utrans-to-condmass}
  \condmass_{\utrans}(\cdot\vert x)
  \eqdef\set{q(\cdot\vert x)\in\simplex_\states}
  {(\forall h\in\allgambles(\states))\ex_{q(\cdot\vert x)}(h)\leq\utrans h(x)},
\end{equation}
and then also
\begin{equation}\label{eq:condmass-back-to-utrans}
  \utrans h(x)
  =\max\set{\ex_{q(\cdot\vert x)}(h)}
  {q(\cdot\vert x)\in\condmass_{\utrans}(\cdot\vert x)}.
\end{equation}
With these credal sets, we can associate a set of transition matrices $\transmats_{\utrans}$:
\begin{equation}\label{eq:utrans-to-transmats}
  \transmats_{\utrans}
  \eqdef\set{\transmat\in\reals^{\states\times\states}}
  {(\forall x\in\states)
    (\exists q(\cdot\vert x)\in\condmass_{\utrans}(\cdot\vert x))
    (\forall y\in\states)
    \transmat_{xy}=q(y\vert x)}.
\end{equation}
In other words, each row $\transmat_{x\cdot}$ of any such transition matrix is formed by the transition probabilities corresponding to some element of $\condmass_{\utrans}(\cdot\vert x)$. 
The elements $\transmat$ of $\transmats_{\utrans}$ are the transition matrices that can be constructed using the one-step information contained in the conditional credal sets $\condmass_{\utrans}(\cdot\vert x)$ and therefore in the (one-step) upper transition operator $\utrans$. 
More generally, the set $\transmats_{\utrans^n}$ contains all $n$-step transition matrices that correspond to the $n$-step upper transition operator $\utrans^n$ (see the Appendix for more details about why we can also consider $\utrans^n$ to be an upper transition operator).

\begin{proposition}\label{prop:hartfiel-and-us}
  Consider a stationary imprecise Markov chain with upper transition operator~$\utrans$ and let $n\in\naturals$. Then
  \begin{enumerate}[(i)]
  \item $\transmats_{\utrans}^n\subseteq\transmats_{\utrans^n}$; 
  \item For all real-valued maps $h$ on $\states$ there is some $\transmat\in\transmats_{\utrans}^n$ such that for all $x\in\states$, $\utrans^nh(x)=(\transmat h)_x$;
  \item\label{prop:hartfiel-and-us:same} For all real-valued maps $h$ on $\states$ and all $x\in\states$,
    \begin{equation}
      \label{eq:hartfiel-and-us}
      \utrans^nh(x)
      =\max\set{(\transmat h)_x}{\transmat\in\transmats_{\utrans}^n}
      \quad\text{ and }\quad
      \ltrans^nh(x)
      =\min\set{(\transmat h)_x}{\transmat\in\transmats_{\utrans}^n}.
    \end{equation}
  \end{enumerate}
\end{proposition}
\noindent
We gather from the following counterexample that for $n>1$, $\transmats_{\utrans}^n$ can be strictly included in $\transmats_{\utrans^n}$.
This shows that the model based on imprecise-probability trees and upper transition operators that we have been using, is more detailed than the Markov set chain model. Nevertheless, as Proposition~\ref{prop:hartfiel-and-us}\eqref{prop:hartfiel-and-us:same} indicates, both models yield very strongly related (if not identical) results as far as the calculation of marginal expectations for $X(n)$ is concerned.

\begin{example}\upshape
  Consider~$\utrans\eqdef(1-\epsilon)\id+\cg\epsilon\max$, where $0\leq\epsilon\leq1$ and~$\id$ is the identity operator, which leaves its argument real-valued map~$h$ unchanged: ${\id h=h}$.
This corresponds to a special case of the contamination models~\eqref{eq:contamination-utrans} discussed in Section~\ref{sec:contamination}. 
For the corresponding $2$-step transition operator, we find that~$\utrans^2=(1-\delta)\id+\cg\delta\max$, with $\delta\eqdef\epsilon(2-\epsilon)$.
\par
Let $\card{\states}=2$, then the sets of corresponding transition matrices are
\begin{equation}
  \transmats_{\utrans}
  =\set{\begin{bmatrix}1-\epsilon_1&\epsilon_1\\\epsilon_2&1-\epsilon_2\end{bmatrix}} 
  {0\leq\epsilon_1,\epsilon_2\leq\epsilon}
  \text{ and }
  \transmats_{\utrans^2}
  =\set{\begin{bmatrix}1-\delta_1&\delta_1\\\delta_2&1-\delta_2\end{bmatrix}} 
  {0\leq\delta_1,\delta_2\leq\delta}.
\end{equation}
We now show that the set $\transmats^2_{\utrans}$ is \emph{strictly} contained in $\transmats_{\utrans^2}$.
Any element of $\transmats_{\utrans}^2$ is given by
\begin{equation}
  \begin{bmatrix}
    1-\epsilon_1&\epsilon_1\\ 
    \epsilon_2&1-\epsilon_2
  \end{bmatrix} 
  \begin{bmatrix}
    1-\epsilon_3&\epsilon_3\\ 
    \epsilon_4 &1-\epsilon_4
  \end{bmatrix} 
  = 
  \begin{bmatrix}
    1-\epsilon_1-\epsilon_3+\epsilon_1\epsilon_3+\epsilon_1\epsilon_4
    &\epsilon_1+\epsilon_3-\epsilon_1\epsilon_3-\epsilon_1\epsilon_4\\ 
    \epsilon_2+\epsilon_4-\epsilon_2\epsilon_4-\epsilon_2\epsilon_3
    &1-\epsilon_2-\epsilon_4+\epsilon_2\epsilon_4+\epsilon_2\epsilon_3
  \end{bmatrix}
\end{equation}
for some $0\leq\epsilon_1,\epsilon_2,\epsilon_3,\epsilon_4\leq\epsilon$, and therefore clearly belongs to $\transmats_{\utrans^2}$. But is is straightforward to check that no choice of $\epsilon_1,\epsilon_2,\epsilon_3,\epsilon_4$ in $[0,\epsilon]$ corresponds to the element of $\transmats_{\utrans^2}$ with $\delta_1=\delta_2=\delta=\epsilon(2-\epsilon)$.~$\blacklozenge$
\end{example}

\Citet{skulj2007} calls a compact set $\transmats$ of transition matrices \newconcept {regular} if there is some $n>0$ such that $T_{xy}>0$ for all $T\in\transmats^n$ and all $x,y\in\states$. He then shows that for such regular~$\transmats$ and for all compact~$\margmass_1$, the corresponding sequence of compact sets~$\margmass_n$ converges in Hausdorff norm to the same compact (and invariant) set $\margmass_{\infty}$.
It follows that for all~$h$ and all compact~$\margmass_1$, the sequence of compact sets $\expects_n(h)$ will converge to the same compact set $\expects_\infty(h)$.
This is a clear generalisation of the classical Perron--Frobenius Theorem~\ref{theo:perron-frobenius-classical}.
But it follows from Proposition~\ref{prop:hartfiel-and-us} that for a given stationary imprecise Markov chain with upper transition operator~$\utrans$, the set~$\transmats_{\utrans}$ is regular in \citeauthor{skulj2007}'s sense if and only if for some $n\in\naturals$, $\ltrans^n\ind{\{y\}}(x)>0$ for all $x,y\in\states$. This is much stronger than even our strongest convergence requirement of regularity (under $\uaccess{}{}$), which only involves the condition $\utrans^n\ind{\{y\}}(x)>0$ for all $x,y\in\states$.
\citeauthor{skulj2007} also proves a convergence result for conservative (too large) approximations of the~$\uex_n$, in the special case of a regular (under~$\uaccess{}{}$) imprecise Markov chain whose upper transition operator is $2$-alternating; see Section~\ref{sec:lower-upper-mass} for further details.
\par
We now turn to \citeauthor{hartfiel1991}'s \citep{hartfiel1991,hartfiel1994,hartfiel1998} results. The strongest general convergence result seems to appear in his book \citep[Sec.~3.2]{hartfiel1998}, where he uses the \newconcept{coefficient of ergodicity} $\tau(\transmat)$ of a transition matrix $\transmat$, defined by
\begin{equation}\label{eq:ergod-coeff}
  \tau(\transmat)
  =\frac{1}{2}\max_{x,y\in\states}\sum_{z\in\states}\abs{\transmat_{xz}-\transmat_{yz}}
  =1-\min_{x,y\in\states}\sum_{z\in\states}\min\{\transmat_{xz},\transmat_{yz}\}.
\end{equation}
A transition matrix is called \newconcept{scrambling} if $\tau(\transmat)<1$.
\citeauthor{hartfiel1991} calls a compact set $\transmats$ of transition matrices \newconcept{product scrambling} if there is some $m\in\naturals$ such that $\tau(\transmat)<1$ for all $\transmat\in\transmats^m$.
He then shows that for such product scrambling $\transmats$ and for all compact $\margmass_1$, the corresponding sequence of compact sets $\margmass_n$ converges in Hausdorff norm to the same compact (and invariant) set $\margmass_{\infty}$.
Again, this is a generalisation of the classical Perron--Frobenius Theorem, and it includes \citeauthor{skulj2007}'s above-mentioned result as a special case.
We believe, however, that this approach, based on the coefficient of ergodicity, has a number of drawbacks that our treatment does not have: the condition seems quite hard to check in practise, and it it is hard to interpret directly. We now also argue that it is too strong, at least from our point of view. 

\begin{proposition}\label{prop:product-scrambling}
  Consider a stationary imprecise Markov chain with upper transition operator~$\utrans$. If~$\transmats_{\utrans}$ is product scrambling, then the chain is regularly absorbing.
\end{proposition}
\noindent
Moreover, as the following counterexample shows,  it is easy to find examples of stationary imprecise Markov chains that are regularly absorbing but for which the corresponding set $\transmats_{\utrans}$ is not product scrambling. Another, perhaps more involved, such counterexample will be presented near the end of Section~\ref{sec:k-out-of-n}.

\begin{example}[Vacuous imprecise Markov chain]\label{ex:vacuous-chain}
  Consider an arbitrary state set $\states$ with at least two elements, and the upper transition operator $\utrans$ defined by $\utrans h=\ind{\states}\max h$ for all real-valued maps $h$ on $\states$. 
The set $\transmats_{\utrans}$ that corresponds to this upper transition operator is the set of \emph{all} transition matrices $\transmats_{\mathrm{all}}$, and consequently $\transmats_{\utrans^n}=\transmats_{\utrans}^n=\transmats_{\mathrm{all}}$ for all $n\in\naturals$ as well. 
\par
Consider the unit transition matrix $\transmat$ defined by $\transmat_{xy}=\delta_{xy}$ [Kronecker delta], so the system remains with probability one in any state $x$ that it is in. 
This $\transmat$ belongs to $\transmats_{\utrans^n}=\transmats_{\mathrm{all}}$ for all $n\in\naturals$, but $\tau(\transmat)=1$, so $\transmats_{\mathrm{all}}$ is not product scrambling. 
\par
But the chain is regularly absorbing! It is even regular (under~$\uaccess{}{}$), in a trivial way: ${\utrans^n\ind{\{y\}}(x)=1}$ for all $n\in\naturals$ and all $x,y\in\states$.
Observe that $\utrans^n=\ind{\states}\max$ and therefore $\uex_\infty=\max$ for all~$\uex_1$.~$\blacklozenge$ 
\end{example}

\section{Examples}\label{sec:examples}
In this section, we indicate how the theory developed in the previous sections can be applied in a number of practical situations.
For each of these, the upper expectations are of some special types that are described in the literature on imprecise probabilities.
We present concrete and explicit examples, as well as a number of simulations.

\subsection{Contamination models}\label{sec:contamination}
Suppose we consider a precise stationary Markov chain, with transition operator $\trans$.
We contaminate it with a vacuous model, i.e., we take a convex mixture with the upper transition operator $\cg\max$ of Example~\ref{ex:vacuous-chain}.
This leads to the upper transition operator $\utrans$, defined by
\begin{equation}\label{eq:contamination-utrans}
  \utrans h=(1-\epsilon)\trans h+\cg\epsilon\max h,
\end{equation}
for all $h\in\allgambles(\states)$, where~$\epsilon$ is some constant in the open real interval $(0,1)$.
The underlying idea is that we consider a specific convex neighbourhood of $\trans$.
Since for all~$x$ in~$\states$, $\min\utrans\ind{\{x\}}={(1-\epsilon)\min\trans\ind{\{x\}}+\epsilon}>0$, this upper transition operator (or the associated imprecise Markov chain) is always
regular (under~$\uaccess{}{}$), regardless of whether~$\trans$ is regular (in the sense of Theorem~\ref{theo:perron-frobenius-classical})!
We infer from Theorem~\ref{theo:convergence} that, whatever the initial upper expectation operator~$\uex_1$ is, the upper expectation operator~$\uex_n$ for the state $X(n)$ at time~$n\in\naturals$ will always converge to the same~$\uex_\infty$.
\par
What is this $\uex_\infty$ is for given~$\trans$ and~$\epsilon$?
For any $n\geq1$,
\begin{align}
  \utrans^nh 
  &= (1-\epsilon)^n\trans^nh
  +\cg\epsilon\smashoperator{\sum_{k=0}^{n-1}}(1-\epsilon)^k\max\trans^kh,\\
  \intertext{and therefore}
  \uex_{n+1}(h)
  &=(1-\epsilon)^n\uex_1(\trans^nh)
  +\epsilon\smashoperator{\sum_{k=0}^{n-1}}(1-\epsilon)^k\max\trans^kh.
  \label{eq:contamination-marginal}
\end{align}
If we now let $n\to\infty$, we see that the limit is indeed independent of the initial upper expectation~$\uex_1$:
\begin{equation}\label{eq:contamination-limit}
  \uex_\infty(h)
  =\epsilon\smashoperator{\sum_{k=0}^{\infty}}
  (1-\epsilon)^k\max\trans^kh.
\end{equation}

\begin{example}[Contaminating a cycle]\upshape
  Consider for instance $\states=\{a,b\}$, and let the precise Markov chain be the cycle with period 2, with transition operator $\trans$ given by $\trans h(a)=h(b)$ and $\trans h(b)=h(a)$.
  Then $\trans^{2n}h=h$ and $\trans^{2n+1}h=\trans h$, and therefore $\max\trans^{2n}h=\max\trans^{2n+1}h=\max h$, whence $\uex_\infty(h)=\max h$.
  So the limit upper expectation is vacuous: we lose all information about the value of $X(n)$ as $n\to\infty$.~$\blacklozenge$
\end{example}

\begin{example}[Contaminating a random walk]\upshape
  Consider a random walk, where $\states=\{a,b\}$ and $\trans h=\cg\frac{h(a)+h(b)}{2}$.
  Then we find that $\uex_\infty(h)=\epsilon\max h+(1-\epsilon)\frac{h(a)+h(b)}{2}$.~$\blacklozenge$
\end{example}

\begin{example}[Another contamination model]\upshape\label{ex:contaminated-transient}
  To illustrate the convergence properties of an imprecise Markov chain, let us look at a simple numerical example.
  Again consider $\states=\{a,b\}$ and let the stationary imprecise Markov chain be defined by an initial credal set $\margmass_1=\set{m\in\simplex_{\{a,b\}}}{0.6\leq m(a)\leq0.9}$, and a contamination model of the type~\eqref{eq:contamination-utrans}, with ${\epsilon=0.1}$, and for which the precise transition operator~$\trans$ is defined by the transition matrix
  \begin{equation*}
    \transmat
    \eqdef
    \begin{bmatrix} 
      q(a\vert a) & q(b\vert a)\\ 
      q(a\vert b) & q(b\vert b)
    \end{bmatrix}
    =
    \begin{bmatrix} 
      0.15 & 0.85\\ 
      0.85 & 0.15
    \end{bmatrix}.
  \end{equation*}
  In~Figure~\ref{fig:contaminated-transient} we have plotted the evolution of $\uex_n(\{a\})$ and $\lex_n(\{a\})$, the upper and lower probability for finding the system in state~$a$ at time $n$, which can be calculated efficiently using Eq.~\eqref{eq:contamination-marginal}.
  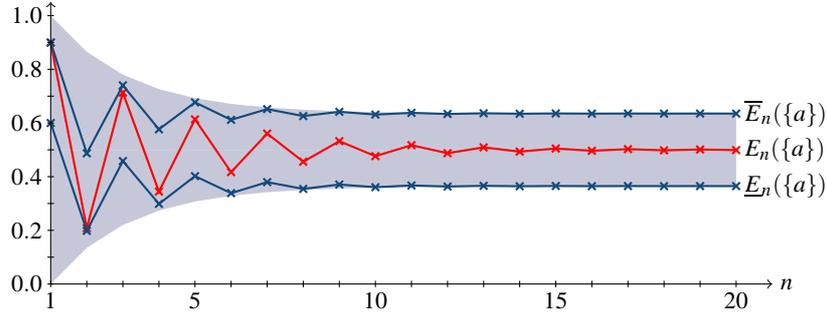
\begin{figure}[ht]
    \centering\small
    \begin{tikzpicture}[baseline,x=1.5em,y=25ex]
      \fill[color=UGentblauw!20] (0,.5) -- plot file {LaV.table} -- (19,.5);
      \fill[color=UGentblauw!20] (0,.5) -- plot file {UaV.table} -- (19,.5);
      \draw[color=red,thick,text=black] plot[mark=x] file {Pa.table} node[right] {$\ex_n(\{a\})$};
      \draw[color=UGentblauw,thick,text=black] plot[mark=x] file {La.table} node[right] {$\lex_n(\{a\})$};
      \draw[color=UGentblauw,thick,text=black] plot[mark=x] file {Ua.table} node[right] {$\uex_n(\{a\})$};
      \draw[->] (0,0) -- (20,0) node[right] {$n$};
      \foreach \k in {0,1,2,3,4,5,6,7,8,9,10,11,12,13,14,15,16,17,18,19} \draw ([yshift=.3ex] \k,0) -- ([yshift=-.3ex] \k,0);
      \foreach \kpos/\k in {0/1,4/5,9/10,14/15,19/20} \path (\kpos,0) node[below] {$\k$};
      \draw[->] (0,0) -- (0,1.05);
      \foreach \val in {0.0,0.1,0.2,0.3,0.4,0.5,0.6,0.7,0.8,0.9,1.0} \draw ([xshift=.3ex] 0,\val) -- ([xshift=-.3ex] 0,\val);
      \foreach \val in {0.0,0.2,0.4,0.6,0.8,1.0} \path (0,\val) node[left] {$\val$};
    \end{tikzpicture}
    \caption{
      The time evolution of (i)~the upper and lower probability of finding the imprecise  Markov chain of Example~\ref{ex:contaminated-transient} in the state~$a$ (outer plot marks and connecting lines); and of (ii)~the probability of finding the classical Markov chain of Example~\ref{ex:contaminated-transient} in the state~$a$ (inner plot marks and connecting lines).
      The filled area denotes the hull of the evolution of this probability, under the contamination model of Example~\ref{ex:contaminated-transient}, for all possible initial mass functions.
    }
    \label{fig:contaminated-transient}
  \end{figure}
  \par
  For comparison, we have also plotted the evolution of~$\ex_n(\{a\})$, the probability for finding the system in state~$a$ at time $n$, for a (precise) Markov chain defined by probability mass functions that lie on the boundaries of the credal sets defining the above imprecise Markov chain; to wit, its initial mass function is given by the row vector~$\dismat_1\eqdef[m_1(a) \;\; m_1(b)]=[0.9 \;\; 0.1]$ and its transition matrix is $\left[\begin{smallmatrix}0.135&0.865\\0.865&0.135\end{smallmatrix}\right]$.
  Here $\ex_\infty(\{a\})=\ex_\infty(\{b\})=0.5$.~$\blacklozenge$
\end{example}

\subsection{Belief function models}\label{sec:belief-models}
The contamination models we have just described are a special case of a more general and quite interesting class of models, based on \citeauthor{shafer1976}'s \citep{shafer1976} notion of a belief function. We can consider a number of subsets $F_j$, $j=1,\dots,n$ of $\states$, and a convex mixture of the vacuous upper expectations relative to these subsets:
\begin{equation}\label{eq:belief-functions}
  \uex(h)=\sum_{j=1}^nm(F_j)\max_{x\in F_j}h(x),
\end{equation}
with $m(F_j)\geq0$ and $\sum_{j=1}^nm(F_j)=1$. In Shafer's terminology, the sets $F_j$ are called \newconcept{focal elements}, and the $m(F_j)$'s the \newconcept{basic probability assignment}.\footnote{Usually, in Shafer's approach, Eq.~\eqref{eq:belief-functions} is only considered for (indicators of) events, and it then defines a so-called \newconcept{plausibility function}, whose conjugate lower probability is a \newconcept{belief function}. Eq.~\eqref{eq:belief-functions} gives the point-wise greatest (most conservative) upper expectation that extends this plausibility function from events to real-valued maps.}
\par
We can now consider imprecise Markov chains where the local models, attached to the non-terminal situations in the tree, are of this type. The general backwards recursion formulae we have given in Section~\ref{sec:sensitivity-analysis} can then be used in combination with the simple formulae of the type~\eqref{eq:belief-functions} for an efficient calculation of all conditional and joint upper and lower expectations in the tree. We leave this implicit however, and move on to another example, which is rather more popular in the literature.

\subsection{Models with lower and upper mass functions}\label{sec:lower-upper-mass}
An intuitive way to introduce imprecise Markov chains \citep{kozine2002,campos2003,skulj2006,hartfiel1998}  goes by way of so-called \newconcept{probability intervals}, studied in a paper by \Citet{campos1994}; see also \citetopt[Section~4.6.1]{walley1991} and \citetopt[Section~2.1]{hartfiel1998}. It consists in specifying lower and upper bounds for mass functions. Let us explain how this is done in the specific context of Markov chains.
\par
For the initial mass function $m_1$, we specify a lower bound $\lmarg_1\colon\states\to\reals$, also called a \newconcept{lower mass function}, and an upper bound $\umarg_1\colon\states\to\reals$, called an \newconcept{upper mass function}.  The credal set $\margmass_1$ attached to the initial situation, which corresponds to these bounds, is then given by
\begin{equation}
  \margmass_1
  \eqdef
  \set{m\in\simplex_\states}
  {(\forall x\in\states)\,\lmarg_1(x)\leq m(x)\leq\umarg_1(x)}.
\end{equation}
\par
Similarly, in each non-terminal situation $\vtuple{x}{k}\in\states^k$, ${k=1,\dots,N-1}$ we have a credal set $\condmass_k(\cdot\vert x_k)$ that is defined in terms of conditional lower and upper mass functions $\lcond_k(\cdot\vert x_k)$ and $\ucond_k(\cdot\vert x_k)$. Here, for instance, $\lcond_k(x_{k+1}\vert x_k)$ gives a lower bound on the transition probability $q_k(x_{k+1}\vert x_k)$ to go from state $X(k)=x_k$ to state $X(k+1)=x_{k+1}$ at time $k$. 
\par
Under some consistency conditions (for more details, see \citep{campos1994}) the upper expectation associated with $\margmass_1$ is then given in all subsets~$A$ of~$\states$ by
\begin{equation}
  \uex_1(A)
  =\min\bigg\{\smashoperator{\sum_{z\in A}}\umarg_1(z),
  1-\smashoperator{\sum_{z\in\states\setminus A}}\lmarg_1(z)\bigg\},
\end{equation}
This $\uex_1$ is \newconcept{$2$-alternating}: $\uex_1(A\cup B)+\uex_1(A\cap B)\leq\uex_1(A)+\uex_1(B)$ for all subsets~$A$ and~$B$ of~$\states$.
This implies (see \citep[Section~3.2.4]{walley1991} and \citep[Theorem~8 and Corollary~17]{cooman2005e}) that for all $h\in\allgambles(\states)$ the upper expectation $\uex_1(h)$ can be found by Choquet integration:
\begin{equation}\label{eq:choquet}
  \uex_1(h)
  =\min h+\smashoperator{\int_{\min h}^{\max h}}
  \uex_1(\set{z\in\states}{h(z)\geq\alpha})\dif\alpha,
\end{equation}
where the integral is a Riemann integral. Similar considerations for the $2$-alternating $\uex_k(\cdot\vert x_k)$ lead to formulae for the upper transition operators $\utrans_k$: for all~$x_k$ in~$\states$,
\begin{align}
  \utrans_k\ind{A}(x_k)
  &=\min\bigg\{\smashoperator{\sum_{z\in A}}\ucond_k(z\vert x_k),
  1-\smashoperator{\sum_{z\in\states\setminus A}}
  \lcond_k(z\vert x_k)\bigg\}\label{eq:choquet2alt}\\
  \utrans_kh(x_k)
  &=\min h+\smashoperator{\int_{\min h}^{\max h}}
  \utrans_k\ind{\set{z\in\states}{h(z)\geq\alpha}}(x_k)
  \dif\alpha.\label{eq:choquet2alttoo}
\end{align}
Using $\uex_1$ and the $\utrans_k$, all (conditional) expectations in the imprecise Markov chain can now be calculated, by applying Theorem~\ref{theo:concatenation} and Corollary~\ref{cor:marginal-concatenation}. 
\par
Rather than using this backwards recursion method, \citet{skulj2006,skulj2007} uses forward propagation, which, reformulated using our notations, amounts to the following. The marginal expectation $\uex_2$ is calculated by $\uex_2=\uex_1\circ\utrans_1$, $\uex_3$ by $\uex_3=\uex_2\circ\utrans_2$, and more generally, $\uex_{n+1}=\uex_n\circ\utrans_n$.
Even though it appears quite natural, this approach has an important drawback, especially in the context of the probability interval approach described above.
In order to calculate, say~$\uex_3(h)$, we first need to find the upper expectation~$\uex_2$, and calculate its value in the map~$\utrans_2h$.
But~$\uex_2$, as the composition of two $2$-alternating models~$\uex_1$ and~$\utrans_1$, is no longer necessarily $2$-alternating, and therefore its value in the map~$\utrans_2h$ cannot generally be calculated from the values it assumes on events, using Choquet integration, as in Eqs.~\eqref{eq:choquet} and~\eqref{eq:choquet2alttoo}.
Indeed, Choquet integration will generally give too large a value for~$\uex_3(h)$, and will therefore lead to conservative approximations.
These are the difficulties that \citeauthor{skulj2006} is faced with in his work \citep{skulj2006,skulj2007}.
\par
They can be circumvented by our backwards recursion approach.
Indeed, in order to find~$\uex_n(h)$, we begin by calculating $h_1\eqdef h$ and $h_{k+1}\eqdef\utrans_kh_k$, $k=1,\dots,n-1$, using Eq.~\eqref{eq:choquet2alttoo}.
Finally, $\uex_n(h)=\uex_1(h_n)$ is calculated using Eq.~\eqref{eq:choquet}.
Our calculations use Choquet integration but are tight, and not conservative approximations, because at all times, the intervening local upper expectations are $2$-alternating.

\begin{example}[Close to a cycle]\label{ex:evosimplex}\upshape
Consider a three-state stationary imprecise Markov model with $\states=\{a,b,c\}$ and with marginal and transition probabilities given by probability intervals.
It follows from Eqs.~\eqref{eq:choquet2alt} and~\eqref{eq:choquet2alttoo} that the upper transition operator~$\utrans$ is fully determined by the lower and upper transition matrices:
\begin{align*}
  \ltransmat\eqdef
  \begin{bmatrix}
    \lcond(a\vert a) & \lcond(b\vert a) & \lcond(c\vert a) \\ 
    \lcond(a\vert b) & \lcond(b\vert b) & \lcond(c\vert b) \\ 
    \lcond(a\vert c) & \lcond(b\vert c) & \lcond(c\vert c)
  \end{bmatrix}
  &=
  \frac{1}{200}
  \begin{bmatrix}
    9 & 9 & 162 \\ 
    144 & 18 & 18 \\ 
    9 & 162 & 9
  \end{bmatrix},\\
  \utransmat\eqdef
  \begin{bmatrix}
    \ucond(a\vert a) & \ucond(b\vert a) & \ucond(c\vert a) \\ 
    \ucond(a\vert b) & \ucond(b\vert b) & \ucond(c\vert b)\\ 
    \ucond(a\vert c) & \ucond(b\vert c) & \ucond(c\vert c)
  \end{bmatrix} 
  &=\frac{1}{200}
  \begin{bmatrix}
    19 & 19 & 172 \\ 
    154 & 28 & 28 \\ 
    19 & 172 & 19
  \end{bmatrix},
\end{align*}
where the numerical values are particular to this example.
We have depicted the credal sets $\condmass(\cdot\vert a)$, $\condmass(\cdot\vert b)$ and $\condmass(\cdot\vert c)$ corresponding to this upper transition operator in Fig.~\ref{fig:near-cyclic-transition}.
\begin{figure}[htb]
  \centering\footnotesize
  \newcommand{\abcsimplex}{
       (0,1,0) node[above] {$c$}
    -- (1,0,0) node[below right] {$b$}
    -- (0,0,1) node[below left] (aunit) {$a$}
    -- cycle
  }
  \begin{tikzpicture}[scale=1.1,z={(-.86603,-.5)},x={(.86603,-.5)},baseline=(aunit)]
    \draw[simplexbackground] \abcsimplex;
    \draw[fill=red] (19/200,162/200,19/200) -- (9/200,172/200,19/200) -- (19/200,172/200,9/200) --  cycle;
    \draw (0.7, 0, 0.7) node {$\condmass(\cdot\vert a)$};
  \end{tikzpicture}
  \hfill
  \begin{tikzpicture}[scale=1.1,z={(-.86603,-.5)},x={(.86603,-.5)},baseline=(aunit)]
    \draw[simplexbackground] \abcsimplex;
    \draw[fill=red] (28/200,18/200,154/200) -- (18/200,28/200,154/200) -- (28/200,28/200,144/200) --  cycle;
    \draw (0.7, 0, 0.7) node {$\condmass(\cdot\vert b)$};
  \end{tikzpicture}
  \hfill
  \begin{tikzpicture}[scale=1.1,z={(-.86603,-.5)},x={(.86603,-.5)},baseline=(aunit)]
    \draw[simplexbackground] \abcsimplex;
    \draw[fill=red] (162/200,19/200,19/200) -- (172/200,19/200,9/200) -- (172/200,9/200,19/200) --  cycle;
    \draw (0.7, 0, 0.7) node {$\condmass(\cdot\vert c)$};
  \end{tikzpicture}
  \hfill
  \begin{tikzpicture}[scale=.7,z={(-.86603,-.5)},x={(.86603,-.5)}]
    \draw (0,1,0) node (c) {$c$} (1,0,0) node (b) {$b$} (0,0,1) node (a) {$a$};
    \draw[->] (c) -- (b);
    \draw[->] (b) -- (a);
    \draw[->] (a) -- (c);
  \end{tikzpicture}
  \caption{The credal sets $\condmass(\cdot\vert a)$, $\condmass(\cdot\vert b)$ and $\condmass(\cdot\vert c)$  in the simplex $\simplex_{\{a,b,c\}}$, corresponding to the upper transition operator $\trans$ in Example~\ref{ex:evosimplex}.}
  \label{fig:near-cyclic-transition}
\end{figure}
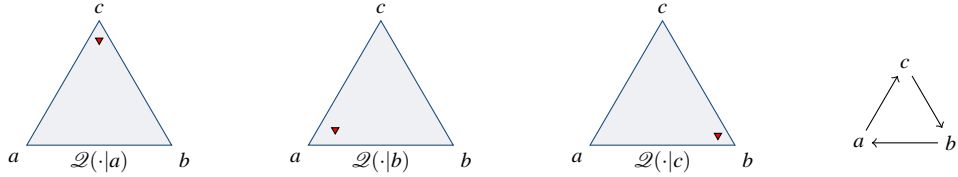

Similarly, the initial upper expectation $\uex_1$ is completely determined by the row vectors $\ldismat_1\eqdef[\lmarg_1(a) \;\; \lmarg_1 (b) \;\; \lmarg_1(c)]$ and $\udismat_1\eqdef[\umarg_1(a) \;\; \umarg_1(b) \;\; \umarg_1(c)]$.
In Figure~\ref{fig:simplex_evolution}, we plot conservative approximations for the credal sets $\margmass_n$ corresponding to the upper expectation operators $\uex_n$.
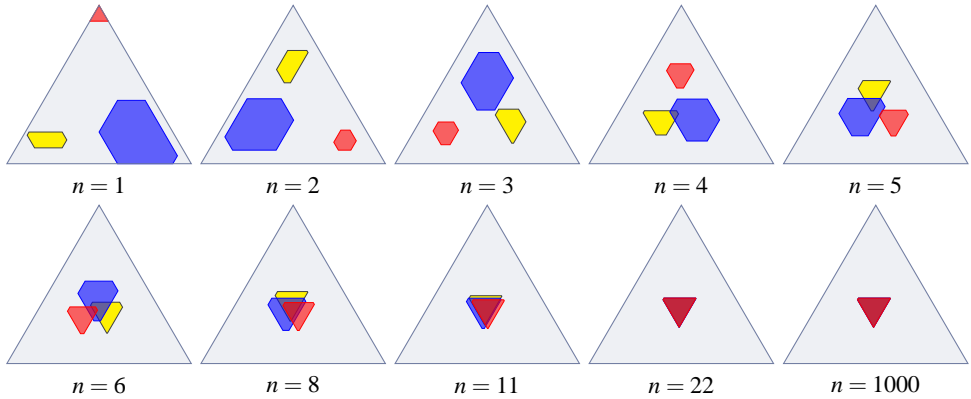
\begin{figure}[htb]
  \centering
  \begin{tikzpicture}[scale=1.4,z={(-.86603,-.5)},x={(.86603,-.5)},baseline=(a)] 
    \fill[simplexbackground] \simplexpath;
    \draw[yellowbackground]     
    (0.1000, 0.1000, 0.8000) -- (0.2500, 0.1000, 0.6500) 
    -- (0.2500, 0.1500, 0.6000) -- (0.2000, 0.2000, 0.6000) 
    -- (0.0200, 0.2000, 0.7800) -- (0.0200, 0.1800, 0.8000) --  cycle;   
    \draw[bluebackground]  
    (0.6000, 0.0000, 0.4000) -- (0.9000, 0.0000, 0.1000) 
    -- (0.9000, 0.0500, 0.0500) -- (0.5500, 0.4000, 0.0500) 
    -- (0.4000, 0.4000, 0.2000) -- (0.4000, 0.2000, 0.4000) --  cycle;   
    \draw[redbackground]  
    (0.0000, 0.9000, 0.1000) -- (0.1000, 0.9000, 0.0000) 
    -- (0.0000, 1.0000, 0.0000) --  cycle;   
    \draw (0.7, 0, 0.7) node {\small$n=1$};   
    \draw[simplexborder] \simplexpath;
  \end{tikzpicture} 
  \hfill 
  \begin{tikzpicture}[scale=1.4,z={(-.86603,-.5)},x={(.86603,-.5)},baseline=(a)] 
    \fill[simplexbackground] \simplexpath;
    \draw[yellowbackground]     
    (0.1125, 0.5945, 0.2930) -- (0.1962, 0.5108, 0.2930) 
    -- (0.2300, 0.5108, 0.2592) -- (0.2300, 0.7016, 0.0684) 
    -- (0.2165, 0.7151, 0.0684) -- (0.1125, 0.7151, 0.1724) --  cycle;   
    \draw[bluebackground]  
    (0.0450, 0.1692, 0.7858) -- (0.1287, 0.0855, 0.7858) 
    -- (0.3650, 0.0855, 0.5495) -- (0.3650, 0.2750, 0.3600) 
    -- (0.2300, 0.4100, 0.3600) -- (0.0450, 0.4100, 0.5450) --  cycle;   
    \draw[redbackground]  
    (0.6525, 0.1355, 0.2120) -- (0.7025, 0.0855, 0.2120) -- 
    (0.7700, 0.0855, 0.1445) -- (0.7700, 0.1445, 0.0855) -- 
    (0.7025, 0.2120, 0.0855) -- (0.6525, 0.2120, 0.1355) --  cycle;   
    \draw (0.7, 0, 0.7) node {\small$n=2$};   
    \draw[simplexborder] \simplexpath;
  \end{tikzpicture} 
  \hfill 
  \begin{tikzpicture}[scale=1.4,z={(-.86603,-.5)},x={(.86603,-.5)},baseline=(a)] 
    \fill[simplexbackground] \simplexpath;
    \draw[yellowbackground]     
    (0.5680, 0.1295, 0.3025) -- (0.5777, 0.1295, 0.2928) 
    -- (0.5777, 0.2644, 0.1579) -- (0.4977, 0.3444, 0.1579) 
    -- (0.3898, 0.3444, 0.2659) -- (0.3898, 0.3078, 0.3025) --  cycle;   
    \draw[bluebackground]  
    (0.1027, 0.5111, 0.3862) -- (0.2750, 0.3389, 0.3862) 
    -- (0.3718, 0.3389, 0.2894) -- (0.3718, 0.5411, 0.0871) 
    -- (0.2107, 0.7022, 0.0871) -- (0.1027, 0.7022, 0.1951) --  cycle;   
    \draw[redbackground]  
    (0.1897, 0.1199, 0.6904) -- (0.2381, 0.1199, 0.6420) 
    -- (0.2381, 0.2116, 0.5503) -- (0.1865, 0.2633, 0.5503) 
    -- (0.1027, 0.2633, 0.6340) -- (0.1027, 0.2069, 0.6904) --  cycle;   
    \draw (0.7, 0, 0.7) node {\small$n=3$};   
    \draw[simplexborder] \simplexpath;
  \end{tikzpicture} 
  \hfill 
  \begin{tikzpicture}[scale=1.4,z={(-.86603,-.5)},x={(.86603,-.5)},baseline=(a)] 
    \fill[simplexbackground] \simplexpath;
    \draw[yellowbackground]     
    (0.2719, 0.1813, 0.5469) -- (0.3275, 0.1813, 0.4912) 
    -- (0.3275, 0.3141, 0.3585) -- (0.3070, 0.3345, 0.3585) 
    -- (0.1324, 0.3345, 0.5331) -- (0.1324, 0.3207, 0.5469) --  cycle;   
    \draw[bluebackground]  
    (0.4578, 0.1433, 0.3989) -- (0.5690, 0.1433, 0.2878) 
    -- (0.5690, 0.2765, 0.1546) -- (0.4375, 0.4080, 0.1546) 
    -- (0.2737, 0.4080, 0.3183) -- (0.2737, 0.3274, 0.3989) --  cycle;   
    \draw[redbackground]  
    (0.2357, 0.4778, 0.2865) -- (0.2727, 0.4778, 0.2495) 
    -- (0.2727, 0.5921, 0.1352) -- (0.2340, 0.6308, 0.1352) 
    -- (0.1260, 0.6308, 0.2433) -- (0.1260, 0.5875, 0.2865) --  cycle;   
    \draw (0.7, 0, 0.7) node {\small$n=4$};   
    \draw[simplexborder] \simplexpath;
  \end{tikzpicture} 
  \hfill 
  \begin{tikzpicture}[scale=1.4,z={(-.86603,-.5)},x={(.86603,-.5)},baseline=(a)] 
    \fill[simplexbackground] \simplexpath;
    \draw[yellowbackground]     
    (0.3080, 0.3335, 0.3585) -- (0.3208, 0.3335, 0.3457) 
    -- (0.3208, 0.5183, 0.1609) -- (0.3141, 0.5250, 0.1609) 
    -- (0.1674, 0.5250, 0.3076) -- (0.1674, 0.4741, 0.3585) --  cycle;   
    \draw[bluebackground]  
    (0.2823, 0.1761, 0.5416) -- (0.3704, 0.1761, 0.4535) 
    -- (0.3704, 0.3603, 0.2693) -- (0.3187, 0.4120, 0.2693) 
    -- (0.1417, 0.4120, 0.4463) -- (0.1417, 0.3167, 0.5416) --  cycle;   
    \draw[redbackground]  
    (0.4956, 0.1752, 0.3292) -- (0.5208, 0.1752, 0.3040) 
    -- (0.5208, 0.3112, 0.1680) -- (0.4941, 0.3379, 0.1680) 
    -- (0.3675, 0.3379, 0.2945) -- (0.3675, 0.3033, 0.3292) --  cycle;   
    \draw (0.7, 0, 0.7) node {\small$n=5$};
    \draw[simplexborder] \simplexpath;
  \end{tikzpicture} 
  \hfill 
  \begin{tikzpicture}[scale=1.4,z={(-.86603,-.5)},x={(.86603,-.5)},baseline=(a)] 
    \fill[simplexbackground] \simplexpath;
    \draw[yellowbackground]     
    (0.4457, 0.1917, 0.3626) -- (0.4494, 0.1917, 0.3589) 
    -- (0.4494, 0.3560, 0.1946) -- (0.4181, 0.3873, 0.1946) 
    -- (0.2701, 0.3873, 0.3426) -- (0.2701, 0.3673, 0.3626) --  cycle; 
    \draw[bluebackground]  
    (0.3371, 0.2696, 0.3933) -- (0.3731, 0.2696, 0.3573) 
    -- (0.3731, 0.4590, 0.1679) -- (0.3118, 0.5203, 0.1679) 
    -- (0.1639, 0.5203, 0.3158) -- (0.1639, 0.4428, 0.3933) --  cycle; 
    \draw[redbackground]  
    (0.3051, 0.1887, 0.5062) -- (0.3231, 0.1887, 0.4882) -- 
    (0.3231, 0.3369, 0.3400) -- (0.3023,       0.3577, 0.3400) 
    -- (0.1633, 0.3577, 0.4790) -- (0.1633, 0.3305, 0.5062) --  cycle; 
    \draw (0.7, 0, 0.7) node {\small$n=6$}; 
    \draw[simplexborder] \simplexpath;
  \end{tikzpicture} 
  \hfill   
  \begin{tikzpicture}[scale=1.4,z={(-.86603,-.5)},x={(.86603,-.5)},baseline=(a)]
    \fill[simplexbackground] \simplexpath;
    \draw[yellowbackground]
    (0.3493, 0.2677, 0.3831) -- (0.3542, 0.2677, 0.3781)
    -- (0.3542, 0.4503, 0.1954) -- (0.3500, 0.4545, 0.1954)
    -- (0.1876, 0.4545, 0.3579) -- (0.1876, 0.4293, 0.3831) --  cycle;
    \draw[bluebackground]
    (0.3395, 0.2094, 0.4511) -- (0.3725, 0.2094, 0.4181)
    -- (0.3725, 0.3918, 0.2357) -- (0.3516, 0.4127, 0.2357)
    -- (0.1779, 0.4127, 0.4095) -- (0.1779, 0.3711, 0.4511) --  cycle;
    \draw[redbackground]
    (0.4188, 0.2088, 0.3724) -- (0.4283, 0.2088, 0.3629)
    -- (0.4283, 0.3734, 0.1983) -- (0.4167, 0.3850, 0.1983)
    -- (0.2618, 0.3850, 0.3532) -- (0.2618, 0.3658, 0.3724) --  cycle;
    \draw (0.7, 0, 0.7) node {\small$n=8$};
    \draw[simplexborder] \simplexpath;
  \end{tikzpicture}
  \hfill 
  \begin{tikzpicture}[scale=1.4,z={(-.86603,-.5)},x={(.86603,-.5)},baseline=(a)] 
    \fill[simplexbackground] \simplexpath;
    \draw[yellowbackground] 
    (0.3645, 0.2432, 0.3922) -- (0.3666, 0.2432, 0.3902) 
    -- (0.3666, 0.4251, 0.2083) -- (0.3633, 0.4284, 0.2083) 
    -- (0.1951, 0.4284, 0.3765) -- (0.1951, 0.4127, 0.3922) --  cycle; 
    \draw[bluebackground] 
    (0.3608, 0.2217, 0.4175) -- (0.3733, 0.2217, 0.4050) 
    -- (0.3733, 0.4034, 0.2233) -- (0.3638, 0.4129, 0.2233) 
    -- (0.1914, 0.4129, 0.3957) -- (0.1914, 0.3911, 0.4175) --  cycle; 
    \draw[redbackground] 
    (0.3903, 0.2213, 0.3884) -- (0.3940, 0.2213, 0.3847) 
    -- (0.3940, 0.3965, 0.2095) -- (0.3880, 0.4025, 0.2095) 
    -- (0.2226, 0.4025, 0.3749) -- (0.2226, 0.3890, 0.3884) --  cycle; 
    \draw (0.7, 0, 0.7) node {\small$n=11$}; 
    \draw[simplexborder] \simplexpath;
  \end{tikzpicture} 
  \hfill 
  \begin{tikzpicture}[scale=1.4,z={(-.86603,-.5)},x={(.86603,-.5)},baseline=(a)] 
    \fill[simplexbackground] \simplexpath;
    \draw[yellowbackground] 
    (0.3732, 0.2287, 0.3981) -- (0.3736, 0.2287, 0.3977) 
    -- (0.3736, 0.4100, 0.2164) -- (0.3708, 0.4128, 0.2164) 
    -- (0.1993, 0.4128, 0.3879) -- (0.1993, 0.4026, 0.3981) --  cycle; 
    \draw[bluebackground] 
    (0.3737, 0.2286, 0.3977) -- (0.3743, 0.2286, 0.3971) 
    -- (0.3743, 0.4099, 0.2158) -- (0.3712, 0.4130, 0.2158) 
    -- (0.1996, 0.4130, 0.3874) -- (0.1996, 0.4026, 0.3977) --  cycle; 
    \draw[redbackground] 
    (0.3731, 0.2295, 0.3974) -- (0.3735, 0.2295, 0.3970) 
    -- (0.3735, 0.4107, 0.2158) -- (0.3707, 0.4136, 0.2158) 
    -- (0.1993, 0.4136, 0.3872) -- (0.1993, 0.4033, 0.3974) --  cycle; 
    \draw (0.7, 0, 0.7) node {\small$n=22$}; 
    \draw[simplexborder] \simplexpath;
  \end{tikzpicture}
  \hfill 
  \begin{tikzpicture}[scale=1.4,z={(-.86603,-.5)},x={(.86603,-.5)},baseline=(a)] 
    \fill[simplexbackground] \simplexpath;
    \draw[yellowbackground] 
    (0.3735, 0.2288, 0.3977) -- (0.3738, 0.2288, 0.3974) 
    -- (0.3738, 0.4102, 0.2160) -- (0.3710, 0.4130, 0.2160) 
    -- (0.1994, 0.4130, 0.3876) -- (0.1994, 0.4028, 0.3977) --  cycle; 
    \draw[bluebackground] 
    (0.3735, 0.2288, 0.3977) -- (0.3738, 0.2288, 0.3974) 
    -- (0.3738, 0.4102, 0.2160) -- (0.3710, 0.4130, 0.2160) 
    -- (0.1994, 0.4130, 0.3876) -- (0.1994, 0.4028, 0.3977) --  cycle; 
    \draw[redbackground] 
    (0.3735, 0.2288, 0.3977) -- (0.3738, 0.2288, 0.3974) 
    -- (0.3738, 0.4102, 0.2160) -- (0.3710, 0.4130, 0.2160) 
    -- (0.1994, 0.4130, 0.3876) -- (0.1994, 0.4028, 0.3977) --  cycle; 
    \draw (0.7, 0, 0.7) node {\small$n=1000$}; 
    \draw[simplexborder] \simplexpath;
  \end{tikzpicture}
  \caption{Evolution in the simplex $\simplex_{\{a,b,c\}}$ of the credal sets $\margmass_n$ for the near-cyclic transition operator from Example~\ref{ex:evosimplex} for three different choices of the initial credal set $\margmass_1$.}
\label{fig:simplex_evolution}
\end{figure}

\noindent Each approximation is based on the constraints that can be found by calculating $\lex_1(\ltrans^{n-1}\ind{\{x\}})$ and $\uex_1(\utrans^{n-1}\ind{\{x\}})$ using the backwards recursion method, for~$x=a,b,c$. 
The~$\margmass_n$ evolve clockwise through the simplex, which is not all that surprising as the lower and upper transition matrices are quite `close' to the precise \emph{cyclic} transition matrix
\begin{equation*}
  \transmat\eqdef
  \begin{bmatrix}
    \cond(a\vert a) & \cond(b\vert a) & \cond(c\vert a) \\ 
    \cond(a\vert b) & \cond(b\vert b) & \cond(c\vert b) \\ 
    \cond(a\vert c) & \cond(b\vert c) & \cond(c\vert c)
  \end{bmatrix} 
  =
  \begin{bmatrix}
    0 & 0 & 1 \\ 1 & 0 & 0 \\ 0 & 1 & 0
  \end{bmatrix},
\end{equation*}
as is also evident from Fig.~\ref{fig:near-cyclic-transition}.
After a while, the $\margmass_n$ converge to a limit that is independent of the initial credal set $\margmass_1$, as can be predicted from the regularity of the upper transition operator.~$\blacklozenge$
\end{example}
\par
A biological application of imprecise Markov models can be found in \citeauthor{dhaenens2007}'s Master's thesis~\citep{dhaenens2007}. He used the sensitivity analysis interpretation of imprecise Markov models to investigate the legitimacy of using \textsc{pam} matrices in amino acid and \textsc{dna} sequence alignments. Roughly speaking, \textsc{pam} (point accepted mutation) matrices describe the chance that one amino acid mutates into another amino acid over a given evolutionary time span. However, the actual value of \textsc{pam} matrix components are based on an estimation using an evolutionary model (i.e., amino acid substitutions are actually counted on the branches of a phylogenetic tree), hence the need to perform a sensitivity analysis. \Citet{dhaenens2007} observed in simulations that the imprecision due to the estimation did not blow up even after a large number of steps; he concluded that using \textsc{pam} matrices over large evolutionary timescales is still reasonable.

\subsection{A \texorpdfstring{$k$}{\itshape k}-out-of-\texorpdfstring{$n$}{\itshape n}:F  system with uncertain reliabilities}\label{sec:k-out-of-n}
Reliability theory is one field where Markov chains are used extensively.
It concerns itself with questions of the type: What is the probability of failure of a system with~$n$ components? 
In the simplest case, where each component is either working or not working, answering this question would involve assessing the failure probabilities of the $2^n$ possible configurations of component states.
However, as shown by \citet{koutras1996}, a great variety of reliability structures can be evaluated quite efficiently using their so-called embedded Markov chain. Amongst these are precisely those systems that fail as soon as any $k$ out of the $n$ components fail, also known as $k$-out-of-$n$:F systems.
\par
For such systems, the embedded Markov chain is constructed as follows. Its state space $\states$ is given by $\{0,1,2,\dots,k\}$, where each number represents the number of components that fail in the system. System failure is therefore represented by the event $\{k\}$, and a fully functioning system by the event $\{0\}$.  \Citet{koutras1996} shows that the failure probability (or unreliability) $F_n$ and the reliability $R_n=1-F_n$ of a  Markov chain embedded system are determined by the expectation form expression: 
\begin{equation}
  F_n 
  \eqdef\ex_{n+1}(\ind{\{k\}})
  =\ex_1(\trans_1\trans_2\ldots\trans_n\ind{\{k\}}),
\end{equation}
where the initial distribution $\ex_1$ represents a system in perfect working condition, so $\ex_1(h)=h(0)$ for all real-valued maps $h$ on $\states$. The transition matrix $\transmat_i$ corresponding to the transition operator $\trans_i$ is fully determined by the reliability $p_i$ of the $i$-th component:
\begin{equation}
  \transmat_i= 
  \begin{bmatrix}
    p_i&1-p_i&0&\dots&0&0\\ 
    0&p_i&1-p_i&\dots&0&0\\ 
    \vdots&\vdots&\vdots&\ddots&\vdots&\vdots\\
    0&0&0&\dots&p_i&1-p_i\\
    0&0&0&\dots&0&1
  \end{bmatrix},
\end{equation}
where $(\transmat_i)_{\ell,m}=\trans_i\ind{\{m\}}(\ell)$ and $\ell,m\in\{0,1,\dots,k\}$.
\par
Precise assessments of the individual reliabilities of the components $p_i$ are often difficult to come by, as for example, they might depend on climatological parameters, age or maybe even on the failure of other (external) components. However, experts might still be able to give conservative bounds on the individual reliabilities $p_i$. In this case, the embedded Markov chain becomes imprecise, but the corresponding bounds on the reliability and unreliability can still be computed by applying our sensitivity analysis formulas derived above:
\begin{equation}
  \overline{F}_n
  =1-\underline{R}_n
  =\ex_1(\utrans_1\utrans_2\ldots\utrans_n\ind{\{k\}})
  \quad\text{ and }\quad
  \underline{F}_n 
  =1-\overline{R}_n
  =\ex_1(\ltrans_1\ltrans_2\ldots\ltrans_n\ind{\{k\}}).
\end{equation}
When this embedded Markov chain is stationary (meaning that the uncertainty models for the reliability of all components are assumed to be the same), the failure probability bounds are simply computed by $\overline{F}_n=\ex_1(\utrans^n\ind{\{k\}})$ and $\underline{F}_n=\ex_1(\ltrans^n\ind{\{k\}})$.
\par
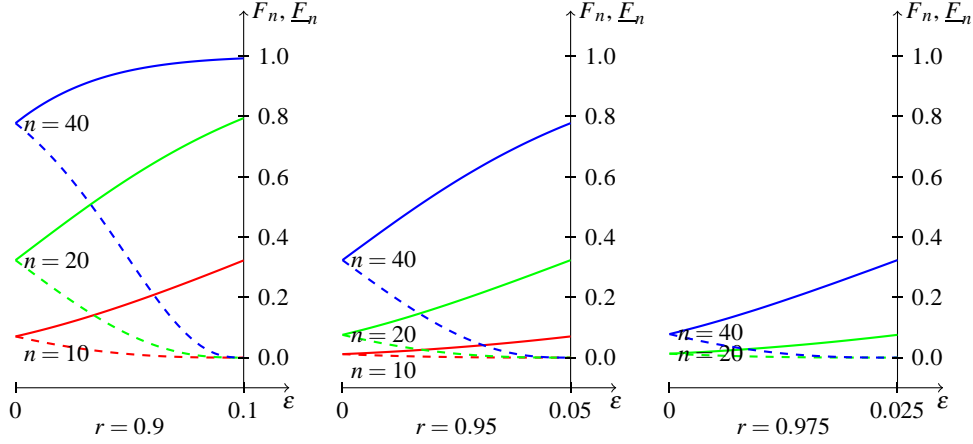
\begin{figure}[hbt]
  \centering\small
  \begin{tikzpicture}
    \foreach \r/\oneminr/\rshift/\rscale in {0.9/0.1/0/1,0.95/0.05/4.3cm/2,0.975/0.025/8.6cm/4} {
      \begin{scope}[xshift=\rshift,yscale=4,xscale=(30*\rscale)]
        \node[below=2ex] at ({.5*(1-\r)},-.1) {$r=\r$};
        \draw[->] (0,-.1) -- ({1.2*(1-\r)},-.1) node[below] {$\varepsilon$};
        \draw[->] ({1-\r},-.1) -- ({1-\r},1.15) node[right] {$\overline{F}_n$, $\underline{F}_n$};
        \foreach \n/\ncolor/\npos in {10/red/below right,20/green/right,40/blue/right} {
          \ifthenelse{\equal{\r\n}{0.97510}}{}{
            \draw[thick,smooth,domain=0:(1-\r),\ncolor] plot[id=u-\r-\n] 
            function{
              1-(\r-x)**\n-\n*(\r-x)**(\n-1)*(1-(\r-x))
              -.5*\n*(\n-1)*(\r-x)**(\n-2)*(1-(\r-x))**2
            };
            \draw[thick,smooth,domain=0:(1-\r),\ncolor,dashed] plot[id=l-\r-\n] 
            function{
              1-(\r+x)**\n-\n*(\r+x)**(\n-1)*(1-(\r+x))  
              -.5*\n*(\n-1)*(\r+x)**(\n-2)*(1-(\r+x))**2
            };
            \node[\npos] at 
            (0,{1-(\r)^\n-\n*(\r)^(\n-1)*(1-\r)-.5*\n*(\n-1)*(\r)^(\n-2)*(1-\r)^2}) {$n=\n$};
          }
        }
        \node[fill,inner sep=0pt,minimum height=1ex,minimum width=.4pt,label=below:$0$] 
        at (0,-.1) {};
        \node[fill,inner sep=0pt,minimum height=1ex,minimum width=.4pt,label=below:$\oneminr$] 
        at ({1-\r},-.1) {};
        \foreach \val in {0.0,0.2,0.4,0.6,0.8,1.0} 
        \node[fill,inner sep=0pt,minimum height=.4pt,minimum width=1ex,label=right:$\val$] 
        at ({1-\r},\val) {};
      \end{scope}
    }
  \end{tikzpicture}
  \caption{
    Upper failure probability ($\overline{F}_n$, full line) and lower failure probability ($\underline{F}_n$, dashed line) for a $3$-out-of-$n$:F system, for different numbers of components~$n$ as a function of the imprecision $\varepsilon\eqdef(\urelty-\lrelty)/2$ of the component reliability, for three different values of $\relty\eqdef(\urelty+\lrelty)/2$.
    As can be expected, the failure bounds widen with increasing imprecision, decrease with increasing reliability (characterised by $\relty$), and increase for a greater number of components~$n$.
  }
  \label{fig:unrel-rel}
\end{figure}

To give a very simple example, let us assume that an expert provides the same range $[\lrelty,\urelty]$ for all component failure probabilities $p_i$, where $0\leq\lrelty\leq\urelty\leq1$. 
This leads to a special case of the models considered in Section~\ref{sec:lower-upper-mass}, and if we apply the formulas derived there, we get, after some manipulations that
\begin{equation}
  \utrans h(\ell)
  =
  \begin{cases}
    \lrelty h(\ell)+(1-\urelty)h(\ell+1)+(\urelty-\lrelty)\max\{h(\ell),h(\ell+1)\}
    &\text{if $\ell=0,1,\dots,k-1$}\\
    h(k)&\text{if $\ell=k$}
  \end{cases}
\end{equation}
for all real-valued maps $h$ on $\states$. 
If $h$ is non-decreasing in the sense that $h(0)\leq h(1)\leq\dots\leq h(k-1)\leq h(k)$, then so is $\utrans h$, and it therefore follows that
\begin{align}
  \overline{F}_n 
  &=
  \begin{bmatrix}
    1&0&\dots&0&0
  \end{bmatrix}
  \begin{bmatrix}
    \lrelty&1-\lrelty&0&\dots&0&0\\
    0&\lrelty&1-\lrelty&\dots&0&0\\
    \vdots&\vdots&\vdots&\ddots&\vdots&\vdots\\
    0&0&0&\dots&\lrelty&1-\lrelty\\
    0&0&0&\dots&0&1
  \end{bmatrix}^n
  \begin{bmatrix}
    0\\0\\\vdots\\0\\1
  \end{bmatrix}\\
  &=\sum_{\ell=k}^n\binom{n}{\ell}\lrelty^{n-\ell}(1-\lrelty)^{\ell}
  =1-\sum_{\ell=0}^{k-1}\binom{n}{\ell}\lrelty^{n-\ell}(1-\lrelty)^{\ell},
\end{align}
and there is a completely similar expression for $\underline{F}_n$ where $\urelty$ is substituted for $\lrelty$.
See Fig.~\ref{fig:unrel-rel} for a graphical illustration of these expressions.
\par
If $0<\lrelty\leq\urelty\leq1$, then this stationary imprecise Markov chain is regularly absorbing with regular top class $\{k\}$ (under~$\uaccess{}{}$), and $\lex_\infty(h)=\uex_\infty(h)=h(k)$ for all real-valued maps $h$ on $\states$. Nevertheless, as soon as $\urelty=1$, \citeauthor{hartfiel1998}'s product scrambling condition is no longer satisfied, as the identity matrix will then belong to all $\transmats_{\utrans^n}$. 
\par
The chain ceases to be regularly absorbing if $\lrelty=0$ and $\urelty=1$, and in that case it is easy to see that $\utrans^{k+n}h(m)=\max_{\ell=m}^kh(\ell)$ for all $n\geq0$ and all real-valued maps $h$ on $\states$, and therefore the limit upper expectation $\uex_\infty$ will depend on the initial upper expectation $\uex_1$. For the particular initial expectation $\ex_1$ we use in this example, we see that $\uex_\infty(h)=\max h$.

\subsection{General models}\label{sec:general-models}
When the (conditional) upper expectation operators that define an imprecise Markov chain do not fall into any of the special cases we discussed and illustrated above, recourse must taken to more general calculation rules.
\par
Let us consider the typical case of a credal set~$\mass$ that is specified by giving, for a finite number of real-valued maps~$f$ collected in the set $\domain\subset\allgambles(\states)$, consistent upper bounds~$U(f)$ on the expectations~$\ex(f)$.
Then the upper expectation for any map $h\in\allgambles(\states)$ can be found by solving the following linear program \citep[see, e.g.,][Section~3.1.3]{walley1991}:
\begin{equation}\label{eq:lin-prog}
\begin{aligned}
  \uex_\mass(h)=\min\bigg[\mu+\smashoperator{\sum_{f\in\domain}}\lambda_fU(f)\bigg]
  \quad\text{subject to }\quad
  &h\leq\mu+\smashoperator{\sum_{f\in\domain}}\lambda_fU(f)\\
  \text{where }\quad 
  &\text{$\lambda_f\geq0$ and $\mu\in\reals$.}
\end{aligned}
\end{equation}
\par
As the number of upper expectations to compute, and thus the number of linear programs to solve, increases, it will eventually become profitable to take a second  (dual) approach.
Any credal set~$\mass$ specified by a finite number of constraints (bounds on expectations) is a convex polytope, i.e., has a finite set~$\ext\mass$ of extreme points.
Vertex enumeration algorithms such as the one by \citet{avis1997} can be used to obtain this set of extreme points from the given set of constraints.
We can then use a practical version of Eq.~\eqref{eq:mass-to-luex} to find the corresponding upper expectations, namely \citep[see][Section~3.1.3]{walley1991}:
\begin{equation}\label{eq:lowens-extver}
  \uex_\mass(h)\eqdef\max\set{\ex_q(h)}{q\in\ext\mass}.
\end{equation}
\par
We can now consider imprecise Markov chains where the local models, attached to the non-terminal situations in the tree, are of this type. The general backwards recursion formulae we have given in Section~\ref{sec:sensitivity-analysis} can then be used in combination with the formulae of the type~\eqref{eq:lin-prog} and~\eqref{eq:lowens-extver} for the calculation of all conditional and joint upper and lower expectations in the tree.

\section{Conclusions}
To conclude, we (i)~reflect on what type of convergence results could be obtained for imprecise Markov chains that are not regularly absorbing, (ii)~we pay attention to the important issue of interpretation of imprecise-probability models, and (iii)~we compare \citeauthor{hartfiel1998}'s approach~\citep{hartfiel1998} to our own regarding their practical applicability to deal with expectation problems.
\par
It is a reasonably weak requirement for a stationary imprecise Markov chain with upper transition operator~$\utrans$ to be regularly absorbing, but we have seen that it is strong enough to guarantee that the upper expectation for the state at time~$n$ converges to a uniquely $\utrans$-invariant upper expectation $\uex_\infty$, regardless of the initial upper expectation $\uex_1$.
\par
Even when an imprecise Markov chain is not regularly absorbing, it is not so hard to see that its upper transition operator~$\utrans$ is still \emph{non-expansive} under the supremum norm given for every ${h\in\allgambles(\states)}$ by $\supnorm{h}\eqdef\max\abs{h}$, as
\begin{equation}
  \supnorm{\utrans g-\utrans h}
  \leq\supnorm{\utrans(g-h)}
  \leq\supnorm{g-h}.
\end{equation}
Moreover, the sequence $\supnorm{\utrans^nh}$ is bounded because $\supnorm{\utrans^nh}\leq\supnorm{h}$.
It then follows from non-linear Perron--Frobenius theory \citep{sine1990,nussbaum1998} that the sequence $\utrans^nh$ has a periodic limit cycle.
More precisely, there is a $\xi_h\in\allgambles(\states)$ such that $\utrans^{p_h}\xi_h=\xi_h$ i.e., $\xi_h$ is a \newconcept{periodic point} of~$\utrans$ with (smallest) \newconcept{period}~$p_h$, and such that $\utrans^{np_h}h\to\xi_h$ (point-wise) as $n\to\infty$.
It would be a very interesting topic for further research to study the nature of the periods and periodic points of upper transition operators.
\par
In our discussions, for instance in Section~\ref{sec:sensitivity-analysis}, we have consistently used the sensitivity analysis interpretation of imprecise-probability models such as upper expectations. Upper and lower expectations can also be given another, so-called \emph{behavioural} interpretation, in terms of some subject's dispositions towards accepting risky transactions.
This is for instance \citeauthor{walley1991}'s \citeyearpar{walley1991} preferred approach.
The results we have derived here remain valid on that alternative interpretation, and the concatenation formulae~\eqref{eq:backpropagation-upper-1} and~\eqref{eq:backpropagation-upper-2} can then be shown to be special cases of so-called \emph{marginal extension} procedure \citep{miranda2006b}, which provides the most conservative coherent (i.e., rational) inferences from the local predictive models $\utrans_k$ to general lower and upper expectations.
In another paper~\Citep{cooman2007d}, we give more details about how to approach a process theory using imprecise probabilities on a behavioural interpretation.
\par
On a related matter: the imprecise Markov chains we are considering here can be seen as special \emph{credal networks} \cite{cozman2000,cozman2005,moral2005b}: the generalisation of Bayesian networks to the case where the local models, associated with the nodes of the network, are credal sets. The corresponding `independence' notion that should then be used for the interpretation of the graphical structure of the network is \citeauthor{walley1991}'s \emph{epistemic irrelevance} \citep[Chapter~9]{walley1991}. Interestingly, \citeauthor{hartfiel1991}'s Markov set-chain approach corresponds to special credal nets where the independence concept involved is a different one: that of \emph{strong independence} \cite{cozman2000}. Nevertheless, both approaches yield the same results if we restrict ourselves to calculating the marginal upper expectations for variables $X(n)$, as we have proved in Proposition~\ref{prop:hartfiel-and-us}. But in any case, for the actual calculation of expectations, the set of transition matrices approach suffers from a combinatorial explosion of computational complexity that can be avoided using our upper transition operator approach.

\section*{Acknowledgements}
The authors wish to thank Damjan \v{S}kulj for inspirational discussion and two anonymous referees for helpful suggestions and pointers to the literature.
\par
This paper presents research results of the Belgian Network DYSCO (Dynamical Systems, Control, and Optimisation), funded by the Interuniversity Attraction Poles Programme, initiated by the Belgian State, Science Policy Office.
The scientific responsibility rests with its authors.


\appendix
\section{Proofs}
In this Appendix, we have gathered proofs for the results in the paper.
\par
Before we go on, it will be useful to discuss and collect a number of properties of the upper transition operators associated with imprecise Markov chains. 
They follow immediately from the corresponding properties~\eqref{eq:uex1}--\eqref{eq:uex7} of upper expectations, so we omit the proof.

\begin{proposition}[Properties of upper transition operators]\label{prop:upper-trans}
  Consider an imprecise Markov chain with a set of states $\states$ and upper transition   operators $\utrans_k$. Then for arbitrary $h$, $h_1$, $h_2$, $h_n$ in $\allgambles(\states)$,   real $\lambda\geq0$ and real $\mu$:   \renewcommand\theenumi{$\utrans$\ensuremath{\arabic{enumi}}}
\begin{enumerate}
\item $\cg\min h\leq\utrans_kh\leq\cg\max h$ (boundedness);\label{eq:utrans1}
\item $\utrans_k(h_1+h_2)\leq\utrans_kh_1+\utrans_kh_2$ (subadditivity); \label{eq:utrans2}
\item $\utrans_k(\lambda h)=\lambda\utrans_kh$ (non-negative homogeneity); \label{eq:utrans3}
\item $\utrans_k(h+\mu\cg)=\utrans_kh+\mu\cg$ (constant additivity); \label{eq:utrans4}
\item if $h_1\leq h_2$ then $\utrans_kh_1\leq\utrans_kh_2$ (monotonicity); \label{eq:utrans5}
\item if $h_n\to h$ point-wise then $\utrans_kh_n\to\utrans_kh$ point-wise (continuity);\label{eq:utrans6}
\item $\utrans_kh\geq-\utrans_k(-h)=\ltrans_kh$ (upper--lower consistency).\label{eq:utrans7}
\end{enumerate}
\end{proposition}
\noindent
Consider any operator $\utrans\colon\allgambles(\states)\to\allgambles(\states)$ that satisfies \eqref{eq:utrans1}--\eqref{eq:utrans3}. 
Then for each  $x\in\states$, the real functional $\uex(\cdot\vert x)$ defined on $\allgambles(\states)$ by $\uex(h\vert x)=\utrans h(x)$ is an upper expectation, because it satisfies \eqref{eq:uex1}--\eqref{eq:uex3}. 
This means that we can consider $\utrans$ as an upper transition operator associated with some imprecise Markov chain. It therefore make sense to call any operator $\utrans$ that satisfies  \eqref{eq:utrans1}--\eqref{eq:utrans3} an \newconcept{upper transition operator}. 
Clearly, if $\utrans_1$, \dots $\utrans_n$ are upper transition operators, then so is their composition $\utrans_1\dots\utrans_n$.
\par
We are now ready to proceed with the proofs of all results in the body of the paper.

\begin{proof}[Proof of Theorem~\ref{theo:concatenation}]
  We first prove by induction that the left-hand sides are dominated by the right-hand sides in   Eqs.~\eqref{eq:backpropagation-upper-1} and~\eqref{eq:backpropagation-upper-2}. To get the   induction process started, we observe that Eq.~\eqref{eq:backpropagation-upper-1} holds   trivially for $n=N-1$. Next, we prove that if the desired inequality in   Eq.~\eqref{eq:backpropagation-upper-1} holds for $n=k+1$, it also holds for $n=k$, where $k$   is any element in $\{1,2,\dots,N-2\}$. Let us fix $\vtuple{x}{k}\in\states^k$, then we have to   prove that
  \begin{equation}
    \uex(f\vert\ntuple{x}{k})
    \leq\uttrans_{k}\uttrans_{k+1}\dots\uttrans_{N-1}f\ftuple{x}{k},
  \end{equation}
  where we can use that, in particular, for all $x_{k+1}\in\states$:
  \begin{equation}\label{eq:starting-point}
    \uex(f\vert\ntuple{x}{k},x_{k+1})
    \leq\uttrans_{k+1}\uttrans_{k+2}\dots\uttrans_{N-1}f(\ntuple{x}{k},x_{k+1}).  
  \end{equation}
  We have fixed $\vtuple{x}{k}$, so we can regard $\uex(f\vert\ntuple{x}{k},\cdot)$ as a   real-valued map on $\states$, depending only on the state $X(k+1)$ at time $k+1$. We denote   this map by $h_{k+1}$.
  \par
  Now consider any compatible probability tree. In particular, let   $q(\cdot\vert\ntuple{x}{k})\in\condmass_k(\cdot\vert x_k)$ be the corresponding local   probability mass function for the uncertainty about the state $X(k+1)$ in the situation   $\vtuple{x}{k}$ we are considering. It follows from the Law of Iterated Expectations that in   this probability tree
  \begin{equation}
    \ex(f\vert\ntuple{x}{k})
    =\ex(\ex(f\vert\ntuple{x}{k},\cdot)\vert\ntuple{x}{k}),
  \end{equation}
  and since $\ex(f\vert\ntuple{x}{k},\cdot)\leq\uex(f\vert\ntuple{x}{k},\cdot)=h_{k+1}$, by definition of the upper expectations in the tree, we may derive from the monotonicity of   expectation operators that $\ex(f\vert\ntuple{x}{k})\leq\ex(h_{k+1}\vert\ntuple{x}{k})$. Now,   $h_{k+1}$ is a function of $X(k+1)$ only, so its conditional expectation   $\ex(h_{k+1}\vert\ntuple{x}{k})$ in situation $\vtuple{x}{k}$ can be calculated using the   local conditional model $q(\cdot\vert\ntuple{x}{k})$ for $X(k+1)$, i.e.,
  \begin{equation}
    \ex(h_{k+1}\vert\ntuple{x}{k})
    =\smashoperator{\sum_{x_{k+1}\in\states}}
    h_{k+1}(x_{k+1})q(x_{k+1}\vert\ntuple{x}{k})
    \leq\uex_k(h_{k+1}\vert x_k),
  \end{equation}
  where the inequality follows from Eq.~\eqref{eq:local-upper}. Hence   $\ex(f\vert\ntuple{x}{k})\leq\uex_k(h_{k+1}\vert x_k)$ and therefore
  \begin{align}
    \uex(f\vert\ntuple{x}{k}) 
    &\leq\uex_k(h_{k+1}\vert x_k)
    =\utrans_kh_{k+1}(x_k)\notag\\
    &\leq\utrans_k\left(\uttrans_{k+1}\uttrans_{k+2}
      \dots\uttrans_{N-1}f(\ntuple{x}{k},\cdot)\right)(x_k)
    =\uttrans_k\uttrans_{k+1}\uttrans_{k+2}
    \dots\uttrans_{N-1}f\ftuple{x}{k},
  \end{align}
  where the first inequality follows from the definition of the upper expectations in the tree,   the first equality follows from Eq.~\eqref{eq:trans-upper}, the second inequality from   Eq.~\eqref{eq:starting-point} and the monotonicity~\eqref{eq:utrans5} of upper transition   operators, and the second equality from Eq.~\eqref{eq:trans-upper-general}.
  \par
  In a completely similar way, but now using the model $\margmass_1$ rather than the model   $\condmass_k(\cdot\vert x_k)$, we can prove that the desired inequalities hold for $n=0$, given   that they hold for $n=1$. So now we know that the left-hand sides are dominated by the   right-hand sides in Eqs.~\eqref{eq:backpropagation-upper-1}   and~\eqref{eq:backpropagation-upper-2}.
  \par
  It remains to prove the converse inequalities. Fix any path in the tree. We denote the successive situations on this path by $\init$, $\vtuple{x}{1}$, $\vtuple{x}{2}$, \dots, $\vtuple{x}{N-1}$, $\vtuple{x}{N}$. First, consider the situation $\vtuple{x}{N-1}$ and the partial map $h_N\eqdef f(\ntuple{x}{N-1},\cdot)$, then we   know, because the credal set $\condmass_{N-1}(\cdot\vert x_{N-1})$ is convex and closed, that   there is some probability mass function in $\condmass_{N-1}(\cdot\vert x_{N-1})$, which we   denote by $\hat{q}(\cdot\vert\ntuple{x}{N-1})$, such that
  \begin{align}
    \smashoperator[r]{\sum_{x_N\in\states}}h_N(x_N)\hat{q}(x_N\vert\ntuple{x}{N-1})
    =\uex_{N-1}(h_N\vert x_{N-1})
    &=\utrans_{N-1}f(\ntuple{x}{N-1},\cdot)(x_{N-1})\notag\\[-2.5ex]
    &=\uttrans_{N-1}f\ftuple{x}{N-1},
  \end{align}
  and therefore
  \begin{equation}\label{eq:first-step}
    \uttrans_{N-1}f\ftuple{x}{N-1}
    =\smashoperator[r]{\sum_{x_N\in\states}}
    f(\ntuple{x}{N-1},x_N)\hat{q}(x_N\vert\ntuple{x}{N-1}).
  \end{equation}
  Next, consider the situation~$\vtuple{x}{N-2}$ and the partial map ${h_{N-1}\eqdef\uttrans_{N-1}f(\ntuple{x}{N-2},\cdot)}$. Again we know, since   $\condmass_{N-2}(\cdot\vert x_{N-2})$ is convex and closed, that there is some probability mass function in $\condmass_{N-2}(\cdot\vert x_{N-2})$, which we denote by $\hat{q}(\cdot\vert\ntuple{x}{N-2})$, such that
  \begin{align}
    \smashoperator[r]{\sum_{x_{N-1}\in\states}}
    h_{N-1}(x_{N-1})\hat{q}(x_{N-1}\vert\ntuple{x}{N-2})
    =\uex_{N-2}(h_{N-1}\vert x_{N-2})
    &=\utrans_{N-2}\left(\uttrans_{N-1}f(\ntuple{x}{N-2},\cdot)\right)
    (x_{N-2})\notag\\[-2.5ex]
    &=\uttrans_{N-2}\uttrans_{N-1}f\ftuple{x}{N-2}
  \end{align}
  and therefore
  \begin{equation}\label{eq:second-step-first}
    \smashoperator[r]{\sum_{x_{N-1}\in\states}}
    \uttrans_{N-1}f(\ntuple{x}{N-2},x_{N-1})     
    \hat{q}(x_{N-1}\vert\ntuple{x}{N-2})
    =\uttrans_{N-2}\uttrans_{N-1}f\ftuple{x}{N-2}.
  \end{equation}
  If we combine Eqs~\eqref{eq:first-step} and~\eqref{eq:second-step-first}, we find that
  \begin{equation}\label{eq:second-step-final}
    \smashoperator[r]{\sum_{\ntuple[N-1]{x}{N}\in\states^2}}
    f(\ntuple{x}{N-2},\ntuple[N-1]{x}{N})\hat{q}(x_{N-1}\vert\ntuple{x}{N-2})
    \hat{q}(x_N\vert\ntuple{x}{N-1})
    =\uttrans_{N-2}\uttrans_{N-1}f\ftuple{x}{N-2}.
  \end{equation}
  We can obviously continue in this manner until we reach the root of the tree. We have then   effectively constructed a compatible probability tree for which the associated conditional and joint expectation operators satisfy for all situations ($n=1,\dots,N-1$)
  \begin{align}
    \uex(f\vert\ntuple{x}{n}) \geq\hat{\ex}(f\vert\ntuple{x}{n}) 
    &{}\eqdef{}\smashoperator{\sum_{\vtuple[n+1]{x}{N}\in\states^{N-n}}}
    f(\ntuple{x}{n},\ntuple[n+1]{x}{N})
    \smashoperator{\prod_{k=n}^{N-1}}\hat{q}(x_{k+1}\vert\ntuple{x}{k})
    =\uttrans_n\uttrans_{n+1}\dots\uttrans_{N-1}f\ftuple{x}{n},\\
    \uex(f)
    \geq\hat{\ex}(f)
    &{}\eqdef{}\smashoperator{\sum_{\vtuple{x}{N}\in\states^{N}}}f\ftuple{x}{N}
    \hat{m}_1(x_1)\smashoperator{\prod_{k=1}^{N-1}}\hat{q}(x_{k+1}\vert\ntuple{x}{k})
    =\uex_1(\uttrans_1\uttrans_2\dots\uttrans_{N-1}f).
  \end{align}
  This tells us that the converse inequalities in Eqs.~\eqref{eq:backpropagation-upper-1}   and~\eqref{eq:backpropagation-upper-2} hold as well.
\end{proof}

\begin{proof}[Proof of Proposition~\ref{prop:markov}]
    We use Eq.~\eqref{eq:backpropagation-upper-1}. It is clear from the     definition~\eqref{eq:trans-upper-general} of the $\uttrans_k$ that if~$f$ is     $\{n,n+1,\dots,N\}$-measurable, then $\uttrans_{N-1}f$ is $\{n,n+1,\dots,N-1\}$-measurable,     and then $\uttrans_{N-2}\uttrans_{N-1}f$ is also $\{n,n+1,\dots,N-2\}$-measurable; so by     continuing the induction, we find $\uttrans_{n+1}\dots\uttrans_{N-1}f$ is     $\{n,n+1\}$-measurable, and finally, $\uttrans_n\dots\uttrans_{N-1}f$ is $\{n\}$-measurable.
\end{proof}

\begin{proof}[Proof of Corollary~\ref{cor:marginal-concatenation}]
  We use Eqs.~\eqref{eq:backpropagation-upper-1} and~\eqref{eq:backpropagation-upper-2} with $f$  defined as follows: $f\ftuple{x}{N}\eqdef h(x_n)$ for all $\vtuple{x}{N}\in\states^N$. Then,   also using~\eqref{eq:utrans3}, the non-negative homogeneity of upper transition operators, we   find after subsequently applying $\uttrans_{N-1}$, \dots, $\uttrans_\ell$ that
  \begin{equation}
  \begin{aligned}
    \uttrans_{N-1}f\ftuple{x}{N-1}
    &=\utrans_{N-1}(h(x_n)\cg)(x_{N-1})=h(x_{n})\\[-.5ex]
    &\;\;\vdots\\[-.5ex]
    \uttrans_{n}\dots\uttrans_{N-1}f\ftuple{x}{n}
    &=\utrans_{n}(h(x_n)\cg)(x_{n})=h(x_{n})\\
    \uttrans_{n-1}\dots\uttrans_{N-1}f\ftuple{x}{n-1}
    &=\utrans_{n-1}h(x_{n-1})\\
    \uttrans_{n-2}\dots\uttrans_{N-1}f\ftuple{x}{n-2}
    &=\utrans_{n-2}\utrans_{n-1}h(x_{n-2})\\[-.5ex]
    &\;\;\vdots\\[-.5ex]
    \uttrans_\ell\dots\uttrans_{N-1}f\ftuple{x}{\ell}
    &=\utrans_\ell\utrans_{\ell+1}\dots\utrans_{n-1}h(x_\ell),
  \end{aligned}
  \end{equation}
and therefore $\uttrans_\ell\dots\uttrans_{N-1}f(\ntuple{x}{\ell-1},\cdot)   =\utrans_\ell\utrans_{\ell+1}\dots\utrans_{n-1}h$. Applying Proposition~\ref{prop:markov} then   leads to the first desired equality. If, for $\ell=1$, we now also apply the upper expectation   $\uex_1$ to both sides of this equality, the proof is complete.
\end{proof}

\begin{proof}[Proof of Proposition~\ref{prop:chapman-kolmogorov}]
  As an example, we prove Eq.~\eqref{eq:CKu}, by applying Eq.~\eqref{eq:backpropagation-upper-1} with its parameters chosen as $f=\ind{\{\vtuple[n+1]{x}{m}\}}$ and $N=m$. We then see that for   any $\vtuple{z}{m-1}\in\states^{m-1}$,
\begin{align}
  \uttrans_{m-1}\ind{\{\vtuple[n+1]{x}{m}\}}\ftuple{z}{m-1}
  &=\utrans_{m-1}\left(\ind{\{\vtuple[n+1]{x}{m-1}\}}
    \ftuple[n+1]{z}{m-1}\ind{\{x_m\}}\right)(z_{m-1})\notag\\
  &=\ind{\{\vtuple[n+1]{x}{m-1}\}}\ftuple[n+1]{z}{m-1}
  \utrans_{m-1}\ind{\{x_m\}}(z_{m-1})\notag\\
  &=\ind{\{\vtuple[n+1]{x}{m-1}\}}\ftuple[n+1]{z}{m-1}
  \utrans_{m-1}\ind{\{x_m\}}(x_{m-1}),
\end{align}
where we have used the non-negative homogeneity~\eqref{eq:utrans3} of upper transition   operators. Therefore   $\uttrans_{m-1}\ind{\{\vtuple[n+1]{x}{m}\}}=\ind{\{\vtuple[n+1]{x}{m-1}\}}\utrans_{m-1}\ind{\{x_m\}}(x_{m-1})$. Consequently, for any $\vtuple{z}{m-2}\in\states^{m-2}$,
\begin{align}
  \uttrans_{m-2}\uttrans_{m-1}\ind{\{\vtuple[n+1]{x}{m}\}}\ftuple{z}{m-2}
  &=\utrans_{m-2}\left(\uttrans_{m-1}\ind{\{\vtuple[n+1]{x}{m}\}}
    (\ntuple{z}{m-2})\right)(z_{m-2})\notag\\
  &\mspace{-40mu}=\utrans_{m-2}\left(\ind{\{\vtuple[n+1]{x}{m-2}\}}\ftuple[n+1]{z}{m-2}
    \ind{\{x_{m-1}\}}\utrans_{m-1}\ind{\{x_m\}}(x_{m-1})\right)(z_{m-2})\notag\\
  &\mspace{-40mu}=\ind{\{\vtuple[n+1]{x}{m-2}\}}\ftuple[n+1]{z}{m-2}
  \utrans_{m-1}\ind{\{x_m\}}(x_{m-1})\utrans_{m-2}\ind{\{x_{m-1}\}}(z_{m-2})\notag\\
  &\mspace{-40mu}=\ind{\{\vtuple[n+1]{x}{m-2}\}}\ftuple[n+1]{z}{m-2}
  \utrans_{m-1}\ind{\{x_m\}}(x_{m-1})\utrans_{m-2}\ind{\{x_{m-1}\}}(x_{m-2}),
\end{align}
again using~\eqref{eq:utrans3}, and therefore
\begin{equation}
  \uttrans_{m-2}\uttrans_{m-1}\ind{\{\vtuple[n+1]{x}{m}\}}
  =\ind{\{\vtuple[n+1]{x}{m-2}\}}\utrans_{m-1}\ind{\{x_m\}}(x_{m-1})
  \utrans_{m-2}\ind{\{x_{m-1}\}}(x_{m-2}).
\end{equation}
Continuing in this fashion eventually leads to Eq.~\eqref{eq:CKu}. 
\end{proof}

\begin{proof}[Proof of Proposition~\ref{prop:topclassregular}]
  Suppose $\maxregstates_{\access[]{}{}}\neq\emptyset$. Consider any maximal state $y$ [there always is at least one, because $\states$ is finite] and any $x\in\maxregstates_{\access[]{}{}}$, then it is clear from the definition of $\maxregstates_{\access[]{}{}}$ that $\access{y}{x}$. Since $y$ is maximal, it follows that also $\access{x}{y}$, and therefore $\commun{x}{y}$. We conclude that $\maxregstates_{\access[]{}{}}$ is included in all maximal communication classes. This means that there is only one such maximal class, and $\maxregstates_{\access[]{}{}}$ is included in this top class. To show that $\maxregstates_{\access[]{}{}}$ is equal to this top class, consider any maximal element $y$ and any $x\in\maxregstates_{\access[]{}{}}$. Then we know that there is some $n\in\naturals$ such that for all $k\geq n$ and all $z\in\states$, $\access[k]{z}{x}$. But we have seen above that $\commun{x}{y}$, so there is some $\ell\geq0$ such that $\access[\ell]{x}{y}$, and therefore $\access[k+\ell]{z}{y}$ for all $z\in\states$. This implies that $y\in\maxregstates_{\access[]{}{}}$, so $\maxregstates_{\access[]{}{}}$ is indeed the top class. We show that it is regular.  For each $x$ in $\maxregstates_{\access[]{}{}}$ there is an $n_x\in\naturals$ such that $\access[k]{y}{x}$ for all $k\geq n_x$ and all $y\in\states$. If we define $n\eqdef\max\set{n_x}{x\in\maxregstates_{\access[]{}{}}}$, then we see that $\access[k]{x}{y}$ for all $k\geq n$ and all $x,y\in\maxregstates_{\access[]{}{}}$, so $\maxregstates_{\access[]{}{}}$ is regular by Proposition~\ref{prop:regularity-characterisation}, and therefore $\access[\cdot]{\cdot}{\cdot}$ is top class regular.
  \par
  Conversely, assume that  $\access[\cdot]{\cdot}{\cdot}$ is top class regular. Consider any state $x$ in the top class, and any $y\in\states$. Then there is some $\ell_y\geq0$ such that $\access[\ell_y]{y}{x}$, and it follows from Proposition~\ref{prop:regularity-characterisation} that there is some $n\in\naturals$ such that $\access[k]{x}{x}$ and therefore $\access[\ell_y+k]{y}{x}$ for all $k\geq n$. So if we let $m\eqdef n+\max\set{\ell_y}{y\in\states}$, then we see that $\access[k]{y}{x}$ for all $k\geq m$ and all $y\in\states$, and therefore $x\in\maxregstates_{\access[]{}{}}$, whence $\maxregstates_{\access[]{}{}}\neq\emptyset$. 
\end{proof}

\begin{proof}[Proof of Proposition~\ref{prop:basic-inequality}]
  Fix $x$, $y$ and $z$ in $\states$. Since $\utp[m]{u}{y}=\utrans^{m}\ind{\{y\}}(u)\ge0$ for all $u\in\states$, we have that 
  \begin{equation}
    \utrans^{m}\ind{\{y\}} 
    =\smashoperator{\sum_{u\in\states}} 
    \utrans^{m}\ind{\{y\}}(u)\ind{\{u\}} 
    \geq\utrans^{m}\ind{\{y\}}(z)\ind{\{z\}}.
  \end{equation}
If we now apply the upper transition operator $\utrans$ $n$ times to both sides of this inequality, and repeatedly invoke its monotonicity~\eqref{eq:utrans5} and non-negative homogeneity~\eqref{eq:utrans3}, we find that $\utrans^{n+m}\ind{\{y\}}\geq \utrans^{m}\ind{\{y\}}(z) \utrans^n\ind{\{z\}}$ and hence indeed $\utrans^{n+m}\ind{\{y\}}(x)\geq\utrans^n\ind{\{z\}}(x)\utrans^{m}\ind{\{y\}}(z)$.
\end{proof}

\begin{proof}[Proof of Proposition~\ref{prop:always-arrive-in-a-state}]
  Fix $x$ in $\states$. Boundedness~\eqref{eq:utrans1} and subadditivity~\eqref{eq:utrans2} guarantee that $0<1\leq\utrans^n\ind{\states}(x)\leq \sum_{y\in\states}\utrans^n\ind{\{y\}}(x)$. So there must be some $y\in\states$ for which $\utp[n]{x}{y}=\utrans^n\ind{\{y\}}(x)>0$.
\end{proof}

The following lemma provides a characterisation for top class regularity (under~$\uaccess{}{}$) that is somewhat simpler than the one implicit in Proposition~\ref{prop:topclassregular}.

\begin{lemma}\label{lem:weakly-topclassregular}
  A stationary imprecise Markov chain is top class regular (under~$\uaccess{}{}$) if and only if
  \begin{equation}\label{eq:1}
    \maxregstates_{\uaccess[]{}{}}
    =\set{x\in\states}
    {(\exists n\in\naturals)(\forall y\in\states)\uaccess[n]{y}{x}}\neq\emptyset.
  \end{equation}
\end{lemma}

\begin{proof}
  Let $\maxregstates_{\uaccess[]{}{}}'\eqdef\set{x\in\states}{(\exists n\in\naturals)(\forall y\in\states)\uaccess[n]{y}{x}}$, then by Proposition~\ref{prop:topclassregular} it suffices to prove that $\maxregstates_{\uaccess[]{}{}}=\maxregstates_{\uaccess[]{}{}}'$. It is clear that $\maxregstates_{\uaccess[]{}{}}\subseteq\maxregstates_{\uaccess[]{}{}}'$, so we concentrate on the converse inequality. Consider any  $x\in\states$ and $n\in\naturals$ such that $\uaccess[n]{y}{x}$ for all $y\in\states$. Then it suffices to prove that also $\uaccess[n+1]{y}{x}$ for all $y\in\states$. Fix $y$, then there is some $z\in\states$ such that $\utp[1]{y}{z}>0$, by Proposition~\ref{prop:always-arrive-in-a-state}. But since we know that for this $z$ also $\utp[n]{z}{x}>0$, we infer from Proposition~\ref{prop:basic-inequality} that indeed $\utp[n+1]{y}{x}\geq\utp[1]{y}{z}\utp[n]{z}{x}>0$.
\end{proof}

Before we come to the upper expectation form of the Perron--Frobenius theorem (Theorem~\ref{theo:convergence}), we first prove the following lemmas.

\begin{lemma}\label{lem:limits}
  Let\/ $\utrans$ be an upper transition operator associated with some stationary imprecise Markov chain, meaning that it satisfies~\eqref{eq:utrans1}--\eqref{eq:utrans7}. Consider any $h\in\allgambles(\states)$. Then the real sequence $\min\utrans^nh$, $n\in\naturals$ is non-decreasing and converges to some limit $l(h)\in\reals$. Similarly, the real sequence $\max\utrans^nh$, $n\in\naturals$ is non-increasing and converges to some limit $L(h)\in\reals$. Of course, $\min h\leq l(h)\leq L(h)\leq\max h$.
\end{lemma}

\begin{proof}
  Fix $h$ in $\allgambles(\states)$ and consider any $n$ in $\naturals_0$. From   $\cg[\states]\min\utrans^nh\leq\utrans^nh\leq\cg[\states]\max\utrans^nh$ [by \eqref{eq:utrans1}] we deduce using~\eqref{eq:utrans5} that $\utrans(\cg[\states]\min\utrans^nh)\leq\utrans^{n+1}h\leq\utrans(\cg[\states]\max\utrans^nh)$, and   therefore, using~\eqref{eq:utrans3} and~\eqref{eq:utrans4}, that   $\cg[\states]\min\utrans^nh\leq\utrans^{n+1}h\leq\cg[\states]\max\utrans^nh$. Consequently,
  \begin{equation}\label{eq:max-min-1}
    \min h
    \leq\min\utrans^nh
    \leq\min\utrans^{n+1}h
    \leq\max\utrans^{n+1}h
    \leq\max\utrans ^nh
    \leq\max h.
  \end{equation}
  This tells us that the real sequence $\max\utrans^nh$ is non-increasing and bounded below (by   $\min h$). It therefore converges to some real number $L(h)$. Similarly, the real sequence   $\min\utrans^nh$ is non-decreasing and bounded above (by $\max h$), and therefore converges to   some real number $l(h)$. That $\min h\leq l(h)\leq L(h)\leq\max h$ follows from the inequalities in Eq.~\eqref{eq:max-min-1} by taking the limit $n\to\infty$.  
\end{proof}

\begin{lemma}\label{lem:subsequence}
    Let\/ $\utrans$ be an upper transition operator associated with some stationary imprecise Markov chain, meaning that it satisfies~\eqref{eq:utrans1}--\eqref{eq:utrans7}. Consider any $h\in\allgambles(\states)$. Then there is some $x_o$ in $\states$ such that for all $n\in\naturals$ there is some $k_n>n$ for which $L(h)\leq\utrans^{k_n}h(x_o)$. Moreover, $\lim_{n\to\infty}\utrans^{k_n}h(x_o)=\limsup_{n\to\infty}\utrans^{n}h(x_o)=L(h)$.
\end{lemma}

\begin{proof}
Suppose, \textit{ex absurdo}, that for any $x\in\states$ there is some $n_x\in\naturals$ such that for all $k>n_x$, $\utrans^{k}h(x)<L(h)$. Since $\states$ is finite, this  implies that there is some $n\eqdef\max\set{n_x}{x\in\states}$ such that for all $k>n$, $\max\utrans^{k}h<L(h)$. This contradicts the conclusion $\max\utrans^nh\searrow L(h)$ obtained in Lemma~\ref{lem:limits}. 
\par
Next, we show that $\lim_{n\to\infty}\utrans^{k_n}h(x_o)=L(h)$. For all $n\in\naturals$, $L(h)\leq\utrans^{k_n}h(x_o)\leq\max\utrans^{k_n}h$, and since the subsequence $\max\utrans^{k_n}h$ converges to the same limit $L(h)$ as the convergent sequence $\max\utrans^{n}h$, we see that the sequence $\utrans^{k_n}h(x_o)$ converges to $L(h)$ as well.
\par
To conclude, we show that $\limsup_{n\to\infty}\utrans^{n}h(x_o)=L(h)$. Since the limit superior of a sequence is the supremum of the limits of all its convergent subsequences, and since moreover we have just proved that $\lim_{n\to\infty}\utrans^{k_n}h(x_o)=L(h)$, we infer that $\limsup_{n\to\infty}\utrans^{n}h(x_o)\geq L(h)$. For the converse inequality: starting from $\utrans^{n}h(x_o)\leq\max\utrans^{n}h$ and taking the limit superior on both sides of the inequality yields
$\limsup_{n\to\infty}\utrans^{n}h(x_o)\leq\limsup_{n\to\infty}\max\utrans^{n}h=L(h)$, where the equality follows from Lemma~\ref{lem:limits}.
\end{proof}

\begin{lemma}\label{lem:onelimit}
  Let\/ $\utrans$ be an upper transition operator associated with some stationary imprecise Markov chain, meaning that it satisfies~\eqref{eq:utrans1}--\eqref{eq:utrans7}. Consider any $h\in\allgambles(\states)$.  If the imprecise Markov chain is regularly absorbing, then $l(h)=L(h)$.
\end{lemma}

\begin{proof}
Since the imprecise Markov chain is in particular top class regular (under~$\uaccess{}{}$), we have by Proposition~\ref{prop:topclassregular} that $\maxregstates_{\uaccess[]{}{}}\neq\emptyset$.
Consider any $x\in\maxregstates_{\uaccess[]{}{}}$, then we first prove that $\lim_{n\to\infty}\utrans^nh(x)=l(h)$. We know from the definition of $\maxregstates_{\uaccess[]{}{}}$ that there is some $n_x\in\naturals$ such that $\min\utrans^{n_x}\ind{\{x\}}>0$. Also, for any $n\geq0$,
\begin{equation}
  0
  \leq\left[\utrans^{n}h(x)-\min\utrans^{n}h\right]\ind{\{x\}} 
  \leq\utrans^{n}h-\min\utrans^{n}h,
\end{equation}
and if we apply $\utrans$ $n_x$ times to all sides of these inequalities, we get 
\begin{equation}
  0 
  \leq\left[\utrans^{n}h(x)-\min\utrans^{n}h\right]\utrans^{n_x}\ind{\{x\}} 
  \leq\utrans^{n+n_x}h-\min\utrans^{n}h,
\end{equation}
after repeated use of~\eqref{eq:utrans5}, \eqref{eq:utrans4} and~\eqref{eq:utrans3}. Taking the minimum of all sides of these inequalities leads to
\begin{equation}
  0 
  \leq 
  \left[\utrans^{n}h(x)-\min\utrans^{n}h\right]
  \min\utrans^{n_x}\ind{\{x\}} 
  \leq\min\utrans^{n+n_x}h-\min\utrans^{n}h.
\end{equation}
If we now let $n\to\infty$, we see that since the term on the right converges to zero [see Lemma~\ref{lem:limits}], so must the middle term. Since  $\min\utrans^{n_x}\ind{\{x\}}>0$, this implies that $\utrans^{n}h(x)-\min\utrans^{n}h$ converges to zero, whence indeed $\lim_{n\to\infty}\utrans^nh(x)=\lim_{n\to\infty}\min\utrans^nh=l(h)$.
\par
As a next step, we infer from Lemma~\ref{lem:subsequence} that there is some $x_o$ in $\states$ and some strictly increasing sequence $k_n$ of natural numbers, such that $L(h)\leq\utrans^{k_n}h(x_o)$ for all $n\in\naturals$, and moreover  $\limsup_{n\to\infty}\utrans^{n}h(x_o)=L(h)$.
\par
There are now two possibilities. The first is that $x_o\in\maxregstates_{\uaccess[]{}{}}$. Then it follows from the discussion above that $\lim_{n\to\infty}\utrans^nh(x_o)=l(h)$. But since we also have that $\lim_{n\to\infty}\utrans^nh(x_o)=\lim_{n\to\infty}\utrans^{k_n}h(x_o)=L(h)$, where the last equality follows from Lemma~\ref{lem:subsequence}, we infer that in this case indeed $l(h)=L(h)$.
\par
The second possibility is that $x_o\notin\maxregstates_{\uaccess[]{}{}}$, but then it follows from the assumption that there is some $n_o\in\naturals$ such that $\ltrans^{n_o}\ind{\maxregstates_{\uaccess[]{}{}}}(x_o)>0$. We have for all $n\in\naturals$ that
\begin{equation}
  0 
  \leq \left[\max\utrans^{n}h-\max_{y\in\maxregstates_{\uaccess[]{}{}}}\utrans^{n}h(y)\right]
  \ind{\maxregstates_{\uaccess[]{}{}}} 
  \leq\max\utrans^{n}h - \utrans^{n}h,
\end{equation}
and if we apply $\ltrans$ $n_o$ times to all sides of these inequalities, we get 
\begin{equation}
  0 
  \leq\left[\max\utrans^{n}h-\max_{y\in\maxregstates_{\uaccess[]{}{}}}\utrans^{n}h(y)\right]
  \ltrans^{n_o}\ind{\maxregstates_{\uaccess[]{}{}}}(x_o) 
  \leq\max\utrans^{n}h-\utrans^{n_o+n}h(x_o),
\end{equation}
after repeated use of~\eqref{eq:utrans5}, \eqref{eq:utrans4}, \eqref{eq:utrans3} and~\eqref{eq:utrans7}, some rearranging, and evaluating in $x_o$. If we now take the limit inferior for $n\to\infty$ of all sides in these inequalities, we find:
\begin{equation}\label{eq:max-min-2}
  0 
  \leq \ltrans^{n_o}\ind{\maxregstates_{\uaccess[]{}{}}}(x_o)
  \liminf_{n\to\infty}\left[\max\utrans^{n}h
    -\max_{y\in\maxregstates_{\uaccess[]{}{}}}\utrans^{n}h(y)\right] 
  \leq\liminf_{n\to\infty}\left[\max\utrans^{n}h-\utrans^{n_o+n}h(x_o)\right].
\end{equation}
Since $\max\utrans^nh\to L(h)$ and $\max_{y\in\maxregstates_{\uaccess[]{}{}}}\utrans^{n}h(y)\to l(h)$ [by the reasoning above, $\utrans^nh(y)\to l(h)$ for all $y\in\maxregstates_{\uaccess[]{}{}}$], we infer that $\liminf_{n\to\infty}\left[\max\utrans^{n}h
-\max_{y\in\maxregstates_{\uaccess[]{}{}}}\utrans^{n}h(y)\right]=L(h)-l(h)$ from the properties of the $\liminf$ operator . It also follows for similar reasons that
\begin{equation}
  \liminf_{n\to\infty}\left[\max\utrans^{n}h-\utrans^{n_o+n}h(x_o)\right]
  =\lim_{n\to\infty}\max\utrans^{n}h-\limsup_{n\to\infty}\utrans^{n_o+n}h(x_o)
  =L(h)-L(h).
\end{equation}
So we infer from Eq.~\eqref{eq:max-min-2} that $\ltrans^{n_o}\ind{\maxregstates_{\uaccess[]{}{}}}(x_o)[L(h)-l(h)]=0$, and therefore that also in this case $l(h)=L(h)$, since by assumption $\ltrans^{n_o}\ind{\maxregstates_{\uaccess[]{}{}}}(x_o)>0$.
\end{proof}

\begin{proof}[Proof of Theorem~\ref{theo:convergence}]
  Since $\cg[\states]\min\utrans^nh\leq\utrans^nh\leq\cg[\states]\max\utrans^nh$, and by   Lemma~\ref{lem:onelimit}, both sequences $\min\utrans^nh$ and $\max\utrans^nh$ converge to the same   real limit, which we denote by $\mu_h$, it follows that $\utrans^nh$ converges (point-wise) to $\cg[\states]\mu_h$: $\lim_{n\to\infty}\utrans^nh=\cg\mu_h$. 
If we use   the continuity of the upper expectation operator $\uex_1$, as well as \eqref{eq:utrans4} and   \eqref{eq:utrans3}, we get
\begin{equation}
  \smashoperator{\lim_{n\to\infty}}\uex_1(\utrans^{n-1}h)
  =\uex_1\left(\smashoperator[r]{\lim_{n\to\infty}}\utrans^{n-1}h\right)
  =\uex_1(\cg\mu_h)=\mu_h,
\end{equation}
and this limit is indeed independent of the choice of $\uex_1$. Hence we find for the limit   that $\uex_\infty(h)=\mu_h$. 
\par
To complete the proof, consider any upper expectation $\uex_1$ on $\allgambles(\states)$ and   any $h$ in $\allgambles(\states)$, then for all $n\in\naturals$,   $\uex_1(\utrans^nh)=\uex_1(\utrans^{n-1}\utrans h)$.  If we let $n\to\infty$ on both sides of   this equality, we find that $\uex_\infty(h)=\uex_\infty(\utrans h)$, showing that   $\uex_\infty$ is indeed $\utrans$-invariant. Now let $\uex_{\mathrm{i}}$ be any $\utrans$-invariant upper   expectation on $\allgambles(\states)$. Then we find for any $h$ in $\allgambles(\states)$, and   for all $n\in\naturals$, that $\uex_{\mathrm{i}}(\utrans^{n-1}h)=\uex_{\mathrm{i}}(h)$, and if we let $n\to\infty$ on both sides of this equality, we find that $\uex_\infty(h)=\uex_{\mathrm{i}}(h)$.
\end{proof}

\begin{proof}[Proof of Proposition~\ref{prop:hartfiel-and-us}]
  We begin with the first statement. It clearly suffices to prove that for any $k\in\naturals$, with obvious notations, $\transmats_{\utrans}\cdot\transmats_{\utrans^k}\subseteq\transmats_{\utrans^{k+1}}$. In other words, consider any $R\in\transmats_{\utrans}$ and any $S\in\transmats_{\utrans^k}$, then we have to show that $T\eqdef RS\in\transmats_{\utrans^{k+1}}$. By Eq.~\eqref{eq:utrans-to-transmats}, $R\in\transmats_{\utrans}$ means that for all $x\in\states$ there is some $r(\cdot\vert x)\in\condmass_{\utrans}(\cdot\vert x)$ such that $R_{xy}=r(y\vert x)$ for all $y\in\states$. Similarly, by Eq.~\eqref{eq:utrans-to-transmats}, $S\in\transmats_{\utrans^k}$ means that for all $y\in\states$ there is some $s(\cdot\vert y)\in\condmass_{\utrans^k}(\cdot\vert y)$ such that $S_{yz}=r(z\vert y)$ for all $z\in\states$. Now for all $x\in\states$ and all $h\in\allgambles(\states)$,
  \begin{align*}
    \utrans^{k+1}h(x)
    &=\utrans(\utrans^kh)(x)\\
    &\geq\ex_{r(\cdot\vert x)}(\utrans^k h)
    =\smashoperator{\sum_{y\in\states}}r(y\vert x)\utrans^kh(y)\\
    &\geq\smashoperator{\sum_{y\in\states}}r(y\vert x)\ex_{s(\cdot\vert y)}(h)
    =\smashoperator{\sum_{y\in\states}}r(y\vert x)
    \smashoperator{\sum_{z\in\states}}s(z\vert y)h(z)
    =\smashoperator{\sum_{y,z\in\states}}R_{xy}S_{yz}h(z)
    =\smashoperator{\sum_{z\in\states}}T_{xz}h(z),
  \end{align*}
where both inequalities follow from Eq.~\eqref{eq:utrans-to-condmass}. If we now consider, for each $x\in\states$, the mass function $q(\cdot\vert x)$ given by $q(z\vert x)\eqdef T_{xz}=\sum_{y\in\states}s(z\vert y)r(y\vert x)$ for all $z\in\states$, then this means that $\utrans^{k+1}h(x)\geq\ex_{q(\cdot\vert x)}(h)$ for all $h\in\allgambles(\states)$, and therefore $q(\cdot\vert x)\in\condmass_{\utrans^{k+1}}(\cdot\vert x)$, for all $x\in\states$, by Eq.~\eqref{eq:utrans-to-condmass}. Hence indeed $T\in\transmats_{\utrans^{k+1}}$, by Eq.~\eqref{eq:utrans-to-transmats}.
\par
On to the second statement. We give a proof by induction. We first show that the statement holds for $n=1$. We know from the definition~\eqref{eq:utrans-to-condmass} of $\condmass_{\utrans}(\cdot\vert x)$  and Eq.~\eqref{eq:condmass-back-to-utrans} that for each $x$ in $\states$ there is some $q(\cdot\vert x)\in\condmass_{\utrans}(\cdot\vert x)$ such that $\utrans h(x)=\sum_{y\in\states}q(y\vert x)h(y)$. Therefore the transition matrix $\transmat$, defined by  $\transmat_{xy}\eqdef q(y\vert x)$ for all $x,y\in\states$, belongs to $\transmats_{\utrans}$ [see Eq.~\eqref{eq:utrans-to-transmats}] and satisfies $\utrans h(x)=\sum_{y\in\states}\transmat_{xy}h(y)=(\transmat h)_x$.
\par
Next, we show that if the statement holds for $n=m$ [the induction hypothesis], it also holds for $n=m+1$, where $m\in\naturals$. Consider the real-valued map $g\eqdef\utrans^mh$, then $\utrans^{m+1}h=\utrans g$. We know from the reasoning above that there is some $\transmat_1\in\transmats_{\utrans}$ such that $\utrans g(x)=(\transmat_1g)_x$ for all $x\in\states$. And the induction hypothesis tells us that there is some $\transmat_2\in\smash[b]{\transmats_{\utrans}^m}$ such that $g(y)=\utrans^mh(y)=(\transmat_2h)_y$ for all $y\in\states$. Hence we find that for all $x\in\states$:
\begin{multline}
  \utrans^{m+1}h(x)
  =\utrans g(x)
  =\sum_{y\in\states}(\transmat_1)_{xy}g(y)\\
  =\sum_{y\in\states}(\transmat_1)_{xy}\sum_{z\in\states}(\transmat_2)_{yz}h(z)
  =\sum_{z\in\states}(\transmat_1\transmat_2)_{xz}h(z)
  =(\transmat_1\transmat_2h)_{x},
\end{multline}
and clearly $\transmat_1\transmat_2\in\transmats_{\utrans}^{m+1}$. This concludes the proof of the second statement.
\par
The third statement is an immediate consequence of the first and second statements.
\end{proof}

Finally, we turn to the proof of proposition~\ref{prop:product-scrambling}.
We first prove an alternative characterisation of the product scrambling property.

\begin{lemma}\label{lem:product-scrambling}
  A set~$\transmats$ of transition matrices is product scrambling if and only if
  \begin{equation}\label{eq:product-scrambling}
    (\exists n\in\naturals)
    (\forall k\geq n)
    (\forall\transmat\in\transmats^k)
    (\forall x,y\in\states)
    (\exists z\in\states)
    \transmat_{xz}>0\wedge\transmat_{yz}>0.
  \end{equation}
\end{lemma}

\begin{proof}
  Recall that $\transmats$ is called product scrambling if
  \begin{equation}
    (\exists n\in\naturals)
    (\forall\transmat\in\transmats^n)
    \tau(\transmat)<1.
  \end{equation}
  Since the coefficient of ergodicity satisfies the submultiplicative property \citep[Section~1.2]{hartfiel1998}:
  \begin{equation}
    \tau(\transmat_1\transmat_2)\leq\tau(\transmat_1)\tau(\transmat_2)
    \text{ for all transition matrices $\transmat_1$ and $\transmat_2$},
  \end{equation}
we see that the product scrambling condition is equivalent to [see also \citep[Lemma~3.2]{hartfiel1998} for a related result]:
  \begin{equation}
    (\exists n\in\naturals)
    (\forall k\geq n)
    (\forall\transmat\in\transmats^k)
    \tau(\transmat)<1.
  \end{equation}
Now use Eq.~\eqref{eq:ergod-coeff}.
\end{proof}

\begin{proof}[Proof of Proposition~\ref{prop:product-scrambling}]
  Assume that $\transmats_{\utrans}$ is product scrambling. We prove that this implies that the corresponding stationary imprecise Markov chain with upper transition operator~$\utrans$ is regularly absorbing: (a)~it is top class regular and (b)~for every~$y$ not in the top class $\maxregstates_{\uaccess{}{}}$, there is some $n\in\naturals$ such that $\ltrans^n\ind{\maxregstates_{\uaccess{}{}}}(y)>0$.
  \par
  We first prove that the Markov chain has a top class under $\uaccess{}{}$.
  It follows from the characterisation~\eqref{eq:product-scrambling} of the product scrambling condition in Lemma~\ref{lem:product-scrambling} that
  \begin{equation}\label{eq:product-scrambling-1}
    (\forall x,y\in\states)
    (\exists z\in\states)
    \uaccess{x}{z}\wedge\uaccess{y}{z},
  \end{equation}
  if we also take into account Proposition~\ref{prop:hartfiel-and-us}.
  For any $x,y\in C$, where $C\subseteq\states$ is the [always non-empty] set of all maximal states, we know that $\uaccess{x}{z}\then\uaccess{z}{x}$ and $\uaccess{y}{z}\then\uaccess{z}{y}$ for all $z\in\states$, so we infer from Eq.~\eqref{eq:product-scrambling-1} that both $\uaccess{x}{y}$ and $\uaccess{y}{x}$, so $x$ and $y$ communicate. This means that the whole of~$C$ forms one single communication class: $C$ is the top class.
  \par
  We now show that this top class $C$ is regular, i.e., consists of a single cyclic subclass, if we recall our discussion of periodicity in Section~\ref{sec:accessibility-abstract}. Let $\period{C}$ be the period of the top class $C$, and consider any $x$ and $y$ in $C$. Using the same reasoning as above, we infer from Eq.~\eqref{eq:product-scrambling} and Proposition~\ref{prop:hartfiel-and-us} that for large enough $k$:
  \begin{equation}
    (\exists z_k\in C)
    \uaccess[k]{x}{z_k}\wedge\uaccess[k]{y}{z_k}
  \end{equation}
  [that $z_k\in C$ follows from the fact that $x$ and $y$ are maximal].
Moreover, Proposition~\ref{prop:class-cycle} tells us that for large enough $\ell$ and $\ell'$, $\steps{z_k}{x}+\ell\period{C}\in\nsteps{z_k}{x}$ and $\steps{z_k}{y}+\ell'\period{C}\in\nsteps{z_k}{y}$, and therefore also  $k+\steps{z_k}{x}+\ell\period{C}\in\nsteps{x}{x}$ and $k+\steps{z_k}{y}+\ell'\period{C}\in\nsteps{y}{y}$. This implies that $\steps{z_k}{x}=\steps{z_k}{y}$, and therefore $\steps{x}{y}=0$: $x$~and~$y$ belong to the same cyclic class. This holds for all $x,y\in C$, so~$C$ consists of only one cyclic class (under~$\uaccess{}{}$). The top class $C$ is in other words aperiodic and therefore regular. This proves (a).
  \par
  To prove (b), assume the stationary imprecise Markov chain is top class regular but not regularly absorbing. We show that the set of transition matrices $\transmats_{\utrans}$ cannot be product scrambling. By Definition~\ref{def:regabs}, we know that there is some $y_0\in\states\setminus\maxregstates_{\uaccess{}{}}$ such that $\ltrans^n\ind{\maxregstates_{\uaccess{}{}}}(y_0)=0$ for all $n\in\naturals$. If we now also invoke Eq.~\eqref{eq:hartfiel-and-us} in Proposition~\ref{prop:hartfiel-and-us}, we see that for all $n\in\naturals$, there is some $\transmat^*_n\in\transmats_{\utrans}^n$ such that:
  \begin{equation}\label{eq:scrambling-leaky1}
    (\forall u\in\maxregstates_{\uaccess{}{}})
    (\transmat^*_n)_{y_0u}=0.
  \end{equation}
Now consider any $x_0$ in the top class $\maxregstates_{\uaccess{}{}}$ [this is possible since by assumption $\maxregstates_{\uaccess{}{}}\neq\emptyset$]. Since $x_0$ cannot communicate with any element outside $\maxregstates_{\uaccess{}{}}$, we infer in particular from Eq.~\eqref{eq:hartfiel-and-us} in Proposition~\ref{prop:hartfiel-and-us} that for all $n\in\naturals$:
\begin{equation}\label{eq:scrambling-leaky2}
  (\forall v\in\states\setminus\maxregstates_{\uaccess{}{}})
  (\transmat^*_n)_{x_0v}=0.
\end{equation}
But Eqs.~\eqref{eq:scrambling-leaky1} and~\eqref{eq:scrambling-leaky2} taken together imply [see Eq.~\eqref{eq:ergod-coeff}] that $\tau(\transmat^*_n)=1$ for all $n\in\naturals$, so the set $\transmats_{\utrans}$ is not product scrambling.
\end{proof}

\end{document}